\documentclass[reqno]{amsart}

\usepackage{amssymb,latexsym}

\def\cal{\mathcal}

\def\Bbb{\mathbb}

\def\Im{\text{\rm Im\,}}
\def\Re{\text{\rm Re\,}}

\def\dbar{\overline{\partial}}

\def\bR{\mathbf{R}}
\def\bC{\mathbf{C}}

\def\cL{\mathcal{L}}

\def\sW{\mathcal{W}}

\def\p{\partial} 
\def\ra{\rangle}
\def\la{\langle}

\def\raw{\longrightarrow}

\def\om{\omega}

\def\bpm{\begin{pmatrix}}

\def\endpf{\hfill $\BoxOpTwo$
\medskip

}

\def\raw{\rightarrow}

\newtheorem{thm}{Theorem}[section]   
\newtheorem{proposition}[thm]{Proposition}

\newtheorem{lemma}[thm]{Lemma}

\newtheorem{defn}[thm]{Definition}

 \def\HollowBox #1#2{{\dimen0=#1 \advance\dimen0 by -#2       
       \dimen1=#1 \advance\dimen1 by #2                       
        \vrule height #1 depth #2 width #2                    
        \vrule height 0pt depth #2 width #1                   
        \llap{\vrule height #1 depth -\dimen0 width \dimen1}%
       \hskip -#2                                             
       \vrule height #1 depth #2 width #2}}                   
 \def\BoxOpTwo{\mathord{\HollowBox{6pt}{.4pt}}\,}             

\font\teneufm=eufm10
\font\seveneufm=eufm7
\font\fiveeufm=eufm5
\newfam\eufmfam
\textfont\eufmfam=\teneufm
\scriptfont\eufmfam=\seveneufm
\scriptscriptfont\eufmfam=\fiveeufm

\newfam\msbfam
\font\tenmsb=msbm10  scaled \magstep1 \textfont\msbfam=\tenmsb
\font\sevenmsb=msbm7 scaled \magstep1  \scriptfont\msbfam=\sevenmsb
\font\fivemsb=msbm5  scaled \magstep1   \scriptscriptfont\msbfam=\fivemsb
\def\Bbb{\fam\msbfam \tenmsb}

\def\ZZ{{\Bbb Z}}

\def\j2{\frac{j+1}{2}}
\def\nb{\nu_\beta}
\def\sgn{{\rm sgn}}
\def\dom{{\rm dom}}

\def\b{\beta}

\def\ve{\varepsilon}
\def\d{\delta}
\def\l{\lambda}

\def\om{\omega}

\def\z{\zeta}

\def\x2p{x + 2\pi}
\def\ex2{e^{-x/2}}

\begin{document}

\title[Analysis and geometry on worm domains]{\large Analysis 
and geometry
on 
worm domains}

\author[S. Krantz]{Steven G. Krantz}
\address{American Institute of Mathematics, 360 Portage Avenue,
Palo Alto, California 94306 USA}
\email{{\tt skrantz@aimath.org}}
\vspace*{.15in}
\author[M. Peloso]{Marco M. Peloso}
\address{Dipartimento di Matematica, Corso Duca degli Abruzzi 24,
  Politecnico di Torino, 10129 
Torino\vspace{-.1in}}
\address{
{\it Current address:}\vspace{-.1in}}
\address{Dipartimento di Matematica, 
Universit\`a Statale di Milano, 
Via C. Saldini 50, 20133 Milano, Italy}
\email{{\tt marco.peloso@mat.unimi.it}}

\thanks{{\bf Thanks:} First 
author supported in part by  
a grant from the Dean of the Graduate College at Washington
University and a grant from the National Science Foundation; second author supported
in part by the 2005 Cofin project {\it Analisi Armonica}.  
The first author 
also thanks the American Institute of Mathematics for its hospitality
and support during 
this work.} 
		
\begin{abstract}
In this primarily expository
paper we study the analysis of the Diederich-Forn\ae ss
worm domain in complex Euclidean space.  We review its importance as a domain
with nontrivial Nebenh\"{u}lle, and as a counterexample to a number of
basic questions in complex geometric analysis.  Then we discuss its
more recent significance 
in the theory of partial differential equations:  the worm is the first 
smoothly bounded, pseudoconvex domain
to exhibit global non-regularity for the $\dbar$-Neumann
problem. We take this opportunity to prove a few new facts.
Next we turn to specific properties of the Bergman kernel for the worm
domain. An asymptotic expansion for this kernel is considered, and
applications to function theory and analysis on the worm are provided. 
\end{abstract}

\maketitle

\thispagestyle{empty}

\section{Introduction}

The concept of ``domain of holomorphy'' is central to the
function theory of several complex variables.  The celebrated
solution of the Levi problem tells us that a connected open set
(a {\it domain}) is a domain of holomorphy if and only if it
is pseudoconvex.  For us, in the present paper, pseudoconvexity
is {\it Levi} pseudoconvexity; this is defined in terms of the positive
semi-definiteness of the Levi form.  This notion requires the boundary
of the domain to be at least $C^2$.  When the boundary is not $C^2$ we
can still define a notion of pseudoconvexity that coincides with the
Levi pseudoconvexity in the $C^2$-case.
When the Levi form is positive
{\it definite} 
then we say that the domain is {\it strongly} pseudoconvex.  
The geometry of pseudoconvex domains
has become an integral part of the study of several complex
variables. (See \cite{Kr1} for basic ideas about
analysis in several complex variables.)	  

Consider a pseudoconvex domain $\Omega \subseteq \bC^n$. 
Any such
domain has an exhaustion 
$U_1 \subset \! \subset U_2 \subset \! \subset U_3
\subset \! \subset \cdots \Omega$ with $\cup_j U_j = \Omega$ by
smoothly bounded, strongly 
pseudoconvex domains. 
This information was
fundamental to the solution of the Levi problem (see \cite{Bers} for
this classical approach), and is an important part of the geometric
foundations of the theory of pseudoconvex domains.

It is natural to ask whether there is a dual result for the exterior
of $\Omega$.  Specifically, given a pseudoconvex domain $\Omega$, are
there smoothly bounded, pseudoconvex domains $W_1 \supset \! \supset
W_2 \supset \! \supset W_3 \supset \! \supset \cdots \supset \!
\supset \cdots \overline{\Omega}$ 
such that $\cap_j W_j = \overline{\Omega}$?   A domain having this
property is said to have a {\it Stein neighborhood basis}.  A domain
failing this property is said to 
have {\it nontrivial Nebenh\"{u}lle}.  

Early on, F. Hartogs in 1906 produced the following
counterexample (which has come to be known as the {\it Hartogs
  triangle}):  
Let $\Omega = \{(z_1, z_2) \in \bC^2: 0 < |z_1| < |z_2| < 1\}$.  

\begin{thm}\label{Hartogs-triangle}  \sl
Any function holomorphic on a neighborhood of $\overline{\Omega}$
actually continues analytically 
to $D^2(0,1) \equiv D \times D$.  Thus $\overline{\Omega}$ {\it
  cannot} have a neighborhood basis of pseudoconvex domains.  Instead
it has a nontrivial Nebenh\"{u}lle.
\end{thm} 
\proof
Let $U$ be a neighborhood of $\overline{\Omega}$.  For
$|z_1| < 1$, the analytic discs 
$$
\zeta \mapsto (z_1, \zeta \cdot |z_1|)
$$
have boundary lying in $U$.  But, for $|z_1|$ sufficiently small, the
entire disc lies in $U$.  Thus a  
standard argument (as in the proof of the Hartogs extension
phenomenon---see \cite{Kr1}), sliding 
the discs for increasing $|z_1|$, shows that a holomorphic function on
$U$ will analytically 
continue to $D(0,1) \times D(0,1)$.  That proves the result. 
\endpf

It was, however, believed for many years that the Hartogs example
worked only because 
the boundary 
of $\Omega$ is not smooth (it is only Lipschitz).  Thus, for over
seventy years, mathematicians 
sought a proof that a smoothly bounded pseudoconvex domain {\it will} have a
Stein neighborhood 
basis.	 It came as quite a surprise in 1977 when Diederich and
Forn\ae ss \cite{DFo1} produced  
a smoothly bounded domain---now known as the {\it worm}---which is 
pseudoconvex and which 
does {\it not} have a Stein neighborhood basis.  In fact the
Diederich-Forn\ae ss example is the following.

\begin{defn}\label{Omega-beta}
{\rm
Let $\sW$ denote the 
domain
$$
\sW = \Big\{(z_1, z_2)\in\bC^2:\,  
\big|z_1 - e^{i\log |z_2|^2} \big|^2 < 1 - \eta(\log|z_2|^2) \Big\}
\, , 
$$
where

\begin{itemize}
\item[{\bf (i)}]  $\eta \geq 0$, $\eta$ is even, $\eta$ is convex;\smallskip
\item[{\bf (ii)}]  $\eta^{-1}(0) = I_\mu \equiv [-\mu,\mu]$; \smallskip
\item[{\bf (iii)}]  there exists a number $a > 0$ such that $\eta(x) > 1$ if
    $|x| > a$; \smallskip
\item[{\bf (iv)}]  $\eta'(x) \ne 0$ if $\eta(x) = 1$.
\end{itemize}
}
\end{defn}

Notice that the slices of $\sW$ for $z_2$ fixed are discs
centered on the unit circle with centers that wind $\mu/\pi$ times
about that circle as $|z_2|$ traverses the range of values for which
$\eta(\log|z_2|^2) < 1$. 

It is worth commenting here on the parameter $\mu$
in the definition of $\sW$.  The number $\mu$ 
in some contexts is selected to
be greater than $\pi/2$.
The number $\nu = \pi/2\mu$
 is half
the reciprocal of the number of times that the centers of the circles that
make up the worm traverse their circular path.  

Many authors use the original choice of parameter $\beta$, where
$\mu=\beta-\pi/2$ (see \cite{Ba3, CheS, KrPe} e.g.). Here, we have
preferred
to  use the notation $\mu$, in accord with the sources \cite{Chr1, Chr2}.

\begin{proposition}\label{pseudoconvexity}  \sl
The domain $\sW$ is
smoothly bounded and pseudoconvex.  Moreover, its boundary
is {\it strongly} 
pseudoconvex except at the boundary points $(0,z_2)$
for $\big|\log |z_2|^2 \big|\le\mu$. 
These points constitute an annulus in
$\p\sW$. 
\end{proposition}

\begin{proposition}\label{nontrivial-Nebenhulle} \sl
The smooth worm domain $\sW$ has nontrivial Nebenh\"{u}lle.
\end{proposition}

\noindent The proofs of these propositions are deferred to Section
\ref{2}.
\medskip
		    
As Diederich and Forn\ae ss \cite{DFo1} showed, the worm provides
a counterexample to a number of interesting questions in the
geometric function theory of several complex variables. As an
instance, the worm gives an example of a smoothly bounded,
pseudoconvex domain which lacks a global plurisubharmonic
defining function. It also provides counterexamples in
holomorphic approximation theory. Clearly the worm showed
considerable potential for a central role in the function
theory of several complex variables. But in point of fact the
subject of the worm lay dormant for nearly fifteen years after
the appearance of \cite{DFo1}. It was the remarkable paper of
Kiselman \cite{Ki} that re-established the importance and
centrality of the worm.

In order to put Kiselman's work into context, we must provide a
digression on the subject of biholomorphic mappings of pseudoconvex 
domains.  In the present discussion, all domains $\Omega$ are 
smoothly bounded.  We are interested 
in one-to-one, onto, invertible mappings (i.e., biholomorphic mappings
or biholomorphisms) of domains
$$
\Phi: \Omega_1 \raw \Omega_2 \, .
$$
Thanks to a classical theorem of Liouville (see \cite{KrPa}), there 
are no conformal mappings, other than trivial ones, in higher
dimensional complex Euclidean space.  Thus biholomorphic mappings
are studied instead.
It is well known that the Riemann mapping theorem fails in several
complex variables (see \cite{Kr1, GKr1, GKr2, IKr}).   It is thus 
a matter of considerable interest to find means to classify domains up
to biholomorphic equivalence.\medskip

Poincar\'{e}'s program for such a classification consisted of
two steps: {\bf (1)} to prove that a biholomorphic mapping of
smoothly bounded pseudoconvex domains extends smoothly to a
diffeomorphism of the closures of the domains and {\bf (2)}
to then calculate biholomorphic differential invariants on the
boundary. His program was stymied for more than sixty years
because the machinery did not exist to tackle step {\bf (1)}
The breakthrough came in 1974 with Fefferman's seminal paper
\cite{Fe}. In it he used remarkable techniques of differential
geometry and partial differential equations to prove that a
biholomorphic mapping of smoothly bounded, {\it strongly}
pseudoconvex domains will extend to a diffeomorphism of the
closures.

Fefferman's proof was quite long and difficult, and left open
the question of {\bf (a)} whether there was a more accessible
and more natural approach to the question and {\bf (b)} whether
there were techniques that could be applied to a more general
class of domains.  Steven Bell \cite{Bel1} as well as Bell and Ewa
Ligocka \cite{BelLi} 
provided a compelling answer.\medskip

Let $\Omega$ be a fixed, bounded domain in $\bC^n$. 
Let $A^2(\Omega)$
be the square integrable holomorphic functions on $\Omega$.   Then
$A^2(\Omega)$ is a closed subspace of $L^2(\Omega)$.  The Hilbert
space projection 
$$
P: L^2(\Omega) \raw A^2(\Omega)
$$
can be represented by an integration formula
$$
Pf(z) = \int_\Omega K(z,\zeta) f(\zeta) \, dV(\zeta) \, .
$$
The kernel $K(z,\zeta) = K_\Omega(z,\zeta)$ is called the {\it Bergman
  kernel}.  It 
is an important biholomorphic invariant.  See \cite{Bers, CheS, Kr1} for
all the basic ideas concerning the Bergman kernel.   

Clearly the Bergman projection $P$ is bounded on $L^2(\Omega)$. Notice
that, if $\Omega$ is smoothly bounded, then $C^\infty(\overline{\Omega})$
is dense in $L^2(\Omega)$. 
In fact more is true:  If $\Omega$ is Levi pseudoconvex and smoothly
bounded, 
then $C^\infty(\overline{\Omega}) \cap \{\hbox{holomorphic functions}\}$ is
dense in $A^2(\Omega)$ (see \cite{Cat3}).

Bell \cite{Bel1} has formulated the notion
of Condition $R$ for the domain $\Omega$. We say that $\Omega$ satisfies
{\it Condition $R$} if $P: C^\infty(\overline{\Omega}) \raw
C^\infty(\overline{\Omega})$.  It is known, thanks to the theory
of the $\dbar$-Neumann problem (see Section \ref{ConditionR}),
that strongly pseudoconvex 
domains satisfy Condition $R$.  
Deep work of Diederich-Forn\ae ss \cite{DFo2}
and Catlin \cite{Cat1},\cite{Cat2} shows that domains with real
analytic 
boundary, and also finite type domains, satisfy Condition $R$.  An
important 
formula of Kohn, which we shall discuss in Section \ref{ConditionR},
relates the $\dbar$-Neumann operator 
to the Bergman projection in a useful way (see also \cite{Kr2}).  The
fundamental result of Bell and Bell/Ligocka is as follows.

\begin{thm}\label{bell-ligocka-thm}  \sl
Let $\Omega_j \subset \bC^n$ be  
smoothly bounded, Levi pseudoconvex domains.  Suppose  that
one of the two domains satisfies Condition $R$.	  If
$\Phi: \Omega_1 \raw \Omega_2$ 
is a biholomorphic mapping then $\Phi$ extends to be a $C^\infty$
diffeomorphism of $\overline{\Omega}_1$ to $\overline{\Omega}_2$.
\end{thm}

This result established the centrality of Condition $R$. The techniques of
proof are so natural and accessible that it seems that Condition $R$ is
certainly the ``right'' approach to questions of boundary regularity of
biholomorphic mappings. Work of Boas/Straube in \cite{BoS2} shows that
Condition $R$ is virtually equivalent to natural regularity conditions on
the $\dbar$-Neumann operator.\footnote{Here the
  $\dbar$-Neumann operator 
$N$ is 
the 
natural right inverse to the $\dbar$-Laplacian
$\BoxOpTwo = \dbar^* \dbar + \dbar \dbar^*$; see Section \ref{ConditionR}.}

For later reference, and for its importance in its own right, we
mention here that the above Theorem \ref{bell-ligocka-thm} can be
``localized''.  To be precise, we say that a given smoothly bounded
domain $\Omega$ satisfies the {\em Local Condition R} at a point
$p_0\in\p\Omega$ if there exists a neighborhood $U$ of $p_0$ such that
$P: C^\infty(\overline{\Omega}) \raw L^2(\Omega)\cap
C^\infty(\overline{\Omega}\cap U)$.  Then Bell's local result is as
follows; see \cite{Bel2}.
\begin{thm}\label{local-biholo}
Let $\Omega_j \subset \bC^n$ be  
smoothly bounded, pseudoconvex domains, $j = 1,2$.  Suppose  that $\Omega_1$
satisfies Local Condition $R$ at $p_0\in\p\Omega_1$.	  If
$\Phi: \Omega_1 \raw \Omega_2$ 
is a biholomorphic mapping then
there exists a neighborhood $U$ of $p_0$ such that
 $\Phi$ extends to be a $C^\infty$
diffeomorphism of $\overline{\Omega}_1\cap U $ onto its image.
\end{thm}
\medskip  

We might mention, as important background information, a result of
David Barrett \cite{Ba1} from 1984.  This considerably predates the
work on which 
the present paper concentrates.  It does not concern the worm, but it
does concern the regularity of the Bergman projection.

\begin{thm}  \sl
There exists a smoothly bounded, non-pseudoconvex  
domain $\Omega \subseteq \bC^2$ on which Condition $R$ fails.
\end{thm}

Although Barrett's result is {\it not} on a pseudoconvex domain, it
provides some insight into the trouble that can be caused by rapidly
varying normals to the boundary.   See \cite{Ba2} for some pioneering
work on this idea.\medskip

As indicated above, it was Kiselman \cite{Ki} who established an important
connection between the worm domain and Condition $R$.  He proved that, for
a certain non-smooth version of the worm (see below), a form of Condition
$R$ fails.  

For $s > 0$, let $H^s(\Omega)$ denote the usual Sobolev space on the
domain $\Omega$ (see, for instance, \cite{Hor1}, \cite{Kr2}).
Building on Kiselman's idea, Barrett \cite{Ba3} used an exhaustion
argument to show that the Bergman projection fails to preserve the
Sobolev spaces of sufficiently high order on the {\it smooth} worm.

\begin{thm}\label{Barrett-thm}  \sl
For $\mu>0$, let $\sW$ be the smooth worm,
defined as in Definition \ref{Omega-beta}, and
let $\nu=\pi/2\mu$. Then 
the Bergman projection $P$ on $\sW$
does not map $H^s(\sW)$ to 
$H^s(\sW)$ when $s \geq \nu$.
\end{thm}

The capstone of results, up until 1996, concerning analysis on the worm domain
is the seminal paper of M. Christ.  Christ finally showed that
Condition $R$ fails on the smooth worm.  Precisely, his result is the
following. 

\begin{thm}\label{Christ-thm}  \sl
Let $\sW$ be the smooth worm.  Then there is a function $f
\in C^\infty(\overline{\sW})$ 
such that its Bergman projection $Pf$ is {\it not} in
$C^\infty(\overline{\sW})$. 
\end{thm}

We note explicitly that the result of this theorem is closely tied
to, indeed is virtually equivalent to, the assertion that
the $\dbar$-Neumann problem is not hypoelliptic
on the smooth worm, \cite{BoS2}.

In the work of Kiselman, Barrett, and Christ, the geometry of the
boundary of the worm plays a fundamental role in the analysis. In
particular, the fact that for large $\mu$ the normal rotates quite
rapidly is fundamental to all of the negative results. It is of interest
to develop a deeper understanding of the geometric analysis of the
worm domain, because it will clearly play a seminal role in future
work in the analysis of several complex variables. 
\medskip

We conclude this discussion of biholomorphic mappings with a
consideration of biholomorphic mappings of the worm. It is at
this time unknown whether a biholomorphic mapping of the
smooth worm $\sW$ to another smoothly bounded, pseudoconvex domain
will extend to a diffeomorphism of the closures. Of course the
worm does {\it not} satisfy Condition $R$, so the obvious
tools for addressing this question are not available. As a
partial result, So-Chin Chen has classified the biholomorphic
{\it self}-maps of the worm $\sW$. His result implies
that all biholomorphic
self-maps of the worm $\sW$ {\it do} extend to
diffeomorphisms of the boundary (see Section 6).

This article is organized as follows:

Section \ref{2} gives particulars of the Diederich-Forn\ae ss
worm. Specifically, 
we prove that the worm is Levi pseudoconvex, and we establish that there
is no global plurisubharmonic defining function. We also examine the
Diederich-Forn\ae ss bounded plurisubharmonic exhaustion function on the
smooth worm.

Section \ref{3} considers non-smooth versions of the worm (these originated
with Kiselman).  We outline some of Kiselman's results.

Section \ref{4} discusses the irregularity of the Bergman projection on
the worm.  In particular, we reproduce some of Kiselman's and Barrett's
analysis.

Section \ref{ConditionR} discusses the failure of Condition $R$ on the
worm domains. 

Section \ref{6} treats the automorphism group of the smooth worm.

Section \ref{7} engages in detailed analysis of the non-smooth worms
$D_\beta$ and $D'_\beta$.   Particularly, we treat $L^p$ boundedness
properties of the Bergman projection, and we study the pathology of
the Bergman kernel on the boundary off the diagonal.

Section 8 treats irregularity properties of the Bergman kernel on worm domains.
\bigskip

\section{The Diederich-Forn\ae ss Worm Domain}\label{2}

We now present the details of the first basic properties of the
Diederich-Forn\ae ss worm domain $\sW$.  Recall that 
$\sW$ is defined in Definition \ref{Omega-beta}.  Some material of this
section can also be found in the excellent monograph
\cite{CheS}. 

We begin by proving Proposition \ref{pseudoconvexity}.
\proof[Proof of Proposition \ref{pseudoconvexity}]
Property {\bf (iii)} of the worm shows immediately
that the worm domain is bounded.  Let
\begin{equation}\label{def-rho}
\rho(z_1, z_2) =  
\Big| z_1 + e^{i\log|z_2|^2} \Big|^2 - 1 + \eta(\log |z_2|^2) \, .
\end{equation}
Then $\rho$ is (potentially) a defining function for $\sW$.
If we can show that $\nabla \rho \ne 0$ at each point of $\p
\sW$ then the implicit function theorem guarantees that $\p
\sW$ is smooth. 

If it happens that $\p \rho/\p z_1 (p) = 0$ at some
boundary point $p = (p_1, p_2)$, then we find that
\begin{equation}\label{*}
\frac{\p \rho}{\p z_1} (p) = \overline{p}_1 + e^{-i\log |p_2|^2} = 0 \,
.  
\end{equation}
Now let us look at $\p
\rho/\p z_2$ at the 
point $p$.  Because of (\ref{*}), the first factor in $\rho$
differentiates to 0 and we find that 
$$
\frac{\p \rho}{\p z_2} (p) = \eta'(\log |p_2|^2) \cdot
\frac{\overline{p}_2}{|p_2|^2} \, . 
$$
Since $\rho(p) = 0$, we have that $\eta(\log|p_2|^2) = 1$.   Hence, by 
property {\bf (iv)}, 
$\eta'(\log|p_2|^2) \ne 0$.  
It follows that $\p \rho/\p z_2(p) \ne 0$.  We conclude
that $\nabla \rho(z) \ne 0$ 
for every boundary point $z$.

For the pseudoconvexity, we write
$$
\rho(z) = |z_1|^2 + 2 \Re \bigl(z_1 e^{-i \log |z_2|^2} \bigr) +
\eta(\log|z_2|^2) \, . 
$$
Multiplying through by $e^{\arg z_2^2}$, we have that locally
$\sW$ is
given by 
$$
|z_1|^2 e^{\arg z_2^2} + 2 \Re \bigl( z_1 e^{-i \log z_2^2} \bigr)
+ \eta(\log |z_2|^2) e^{\arg z_2^2} < 0 \, . 
$$
The function $e^{-i\log z_2^2}$ is locally  well defined and 
holomorphic, and its modulus is $e^{\arg z_2^2}$.
Thus the first two terms are plurisubharmonic. 
Therefore we must check that the last term is plurisubharmonic.  Since
it only depends on $z_2$, we merely have to calculate its
Laplacian.  We have, arguing as before, that
$$
\Delta \Bigl( \eta (\log |z_2|^2) e^{\arg z_2^2} \Bigr)  = 
   \Bigl( \Delta \eta(\log |z_2|^2) \Bigr) e^{\arg z_2^2} +
    \eta(\log |z_2|^2) \Delta e^{\arg z_2^2} \ge 0 \, .
$$

Because $\eta$ is convex and nonnegative (property {\bf (i)}), the
nonnegativity 
of this last expression follows.  This shows that $\sW$ is
smoothly bounded and pseudoconvex.

In order to describe the locus
of weakly pseudoconvex points, we consider again the
local defining function 
$$
\rho(z_1, z_2) = |z_1|^2 e^{\arg z_2^2} + 2 \Re \bigl( z_1 e^{-i \log z_2^2} \bigr)
+ \eta(\log |z_2|^2) e^{\arg z_2^2}  \, . 
$$
This function is strictly plurisubharmonic at all points $(z_1,z_2)$
with $z_1\neq0$ because of the first two terms, or  
where $\big|\log |z_2|^2\big|>\mu$, because of the last term.
Thus consider the annulus $\mathcal A\subset\p\sW$ given by
\begin{equation}\label{A}
\mathcal A=\bigl\{ (z_1,z_2)\in\p\sW\, :
\, z_1=0\ \text{and}\ 
\big|\log |z_2|^2\big|\le \mu\,
\bigr\}\ .
\end{equation}
A direct calculation shows that the complex Hessian for $\rho$ at a
point $z\in \mathcal A$ acting on $v=(v_1,v_2)\in\bC^2$ is given by 
$$
|v_1|^2 +
2\Re \Bigl(v_1\bar v_2 \frac{e^{i\log |z_2|^2}}{z_2}\Bigr) \ .
$$
By pseudoconvexity, such an expression must be non-negative for all
complex tangential vectors $v$ at $z$.  But such vectors are of the
form $v=(0,v_2)$, so that the Levi form $\cL_{\rho}\equiv0$ on $A$. 
This proves the result
\endpf

It is appropriate now to give the proof of Diederich and Forn\ae ss
that 
the worm has nontrivial Nebenh\"{u}lle.  What is of interest here, and
what distinguishes the worm from the older example of the Hartogs
triangle, is that the worm is a bounded, pseudoconvex domain with {\it
  smooth   boundary}. \medskip

We now show that 
$\overline{\sW}$ does not have a Stein neighborhood basis.

\proof[Proof of Proposition \ref{nontrivial-Nebenhulle}]
What we actually show is that if $U$ is {\it any} neighborhood
of $\overline{\sW}$ then $U$ will contain
$$
K = \{(0, z_2): -\pi \leq \log |z_2|^2 \leq \pi\} \cup
    \{(z_1, z_2): \log |z_2|^2 = \pi \ \hbox{or} \  -\pi \
    \hbox{and} \ |z_1 - 1| < 1\} \, .
$$
In fact this assertion is immediate by inspection.

By the usual Hartogs extension phenomenon argument, it then follows 
immediately that if $U$ is pseudoconvex then $U$ must contain
$$
\widehat{K} = \bigl\{(0, z_2): -\pi \leq \log |z_2|^2 \leq \pi \
    \hbox{and} \ |z_1 - 1| < 1 \bigr\} \, .
$$
Thus there can be no Stein neighborhood basis.
\endpf

We now turn to a few properties of the smooth worm $\sW$ connected
with potential theory.
The significance of the next result stems from the paper
\cite{Boas-Straube}.
In that paper, Boas and Straube established the following. 
\begin{thm}\label{psh-cond-R}
Let $\Omega$ be 
a smoothly bounded
pseudoconvex domain that admits a defining function that
is plurisubharmonic on the boundary. Then, for every $s>0$,
$$
P:H^s(\Omega)\rightarrow H^s(\Omega)
$$
is bounded.  In particular, 
$\Omega$ satisfies Condition $R$. 
\end{thm}

For sake of completeness we mention here that, if the Bergman
projection $P$ on a domain $\Omega$ is such that
$P:C^\infty(\overline{\Omega}) \raw  C^\infty(\overline{\Omega})$ is
bounded (i.e.
$\Omega$ satisfies Condition $R$) $P$ is said to be {\em
  regular}, while if $P:H^s(\Omega)\rightarrow H^s(\Omega)$ for every
$s>0$ (and hence $\Omega$ satisfies Condition $R$ {\em a fortiori})
$P$ is said to be {\em exactly  regular}.

Thanks to the result of Christ \cite{Chr1}, we now know that $\sW$ 
does not satisfy Condition $R$, hence {\em a fortiori} it cannot
admit a defining function which is plurisubharmonic on the boundary.
However, it is simpler to give a direct proof of this fact.

\begin{proposition}\label{lack-of-psh-def-funct}  \sl
There exists no defining function $\widetilde \rho$ for $\sW$ that
is plurisubharmonic on the entire boundary.
\end{proposition}
\proof
Suppose that such a defining $\widetilde \rho$ exist.  Then, there
exists a smooth positive function $h$ such that $\widetilde
\rho=h\rho$.   
A direct calculation shows that the complex Hessian for
$\widetilde\rho$ at a 
point $z\in \mathcal A$ acting on $v=(v_1,v_2)\in\bC^2$ is given by 
$$
\mathcal L_{\widetilde \rho} (z;(v_1,v_2))= 
2\Re \Bigl[\bar v_1 v_2 \bigl(  \frac{ih}{z_2} + \p_{z_2} h\bigr)
e^{i\log |z_2|^2} \Bigr]
+ \Bigl[ h+2\Re \bigl( \p_{z_1}h \cdot e^{i\log|z_2|^2} \bigr) \Bigr]
 |v_1|^2 \ .
$$
Since this expression is assumed to be always non-negative, we must
have
$$
\bigl(  \frac{ih}{z_2} + \p_{z_2} h\bigr)
e^{i\log z_2^2}  = \overline{\p_{\bar z_2}
\bigl(he^{-i\log|z_2|^2}\bigr)}
\equiv 0\ ,
$$
on $\mathcal A$.  Therefore, 
the function $g(z_2)=h(0,z_2)e^{-i\log|z_2|^2}$ is a holomorphic
function on $\mathcal A$.  Hence $g(z_2)e^{i\log|z_2|^2}=h(0,z_2)
e^{2\arg z_2}$ is locally a holomorphic function.  Thus it must be
locally a constant, hence a constant $c$ on all of $\mathcal A$.

Therefore, on $\mathcal A$,
$$
h(0,z_2)= c
e^{-2\arg z_2}
$$
which is impossible.  This proves the result.
\endpf

We conclude this section with another important result about the
Diederich-Forn\ae ss worm domain $\sW$.  This result is part of 
potential theory, and is related to the negative result
Proposition \ref{lack-of-psh-def-funct}.  
In what follows, we say that $\lambda$ is a 
{\it bounded plurisubharmonic exhaustion function} for a domain
$\Omega$ if 
\begin{enumerate}
\item[{\bf (a)}]  $\lambda$ is continuous on $\overline{\Omega}$;
\item[{\bf (b)}]  $\lambda$ is strictly plurisubharmonic on $\Omega$;
\item[{\bf (c)}]  $\lambda = 0$ on $\p \Omega$;
\item[{\bf (d)}]  $\lambda < 0$ on $\Omega$;
\item[{\bf (e)}]  For any $c < 0$, the set $\Omega_c = \{z \in \Omega:
\lambda(z) < c\}$ is relatively compact in $\Omega$.
\end{enumerate}
A bounded plurisubharmonic exhaustion function carries important
geometric information 
about the domain $\Omega$.

Now Diederich-Forn\ae ss have proved the following \cite{DFo2}.

\begin{thm}\label{existence-bdd-psh-exhaustion-fnct}  \sl
Let $\Omega$ be any smoothly bounded pseudoconvex domain,
$\Omega=\{z\in\bC\, :\ \varrho(z)<0\}$.  Then there exists $\d$,
$0<\d\le1$, 
and a defining function $\widetilde \varrho$ for $\Omega$ such that
$-(-\widetilde\varrho)^\d$ is a bounded strictly plurisubharmonic
exhaustion function for $\Omega$.
\end{thm}

The importance of this result in the setting of the regularity of
the Bergman projection appears in the following related result, proved
by
Berndtsson-Charpentier \cite{BeCh} and Kohn \cite{Kohn-exp},
respectively. 
\begin{thm} 
Let $\Omega$ be a smoothly bounded pseudoconvex domain and let 
$P$ denotes its Bergman projection.  Let $\widetilde{\rho}$ be
a smooth defining function for $\Omega$ such that 
$-(-\widetilde{\rho})^\delta$ is
strictly plurisubharmonic.  Then  there exists
$s_0=s_0(\Omega,\delta)$ such that 
$$
P:H^s(\Omega)\rightarrow H^s(\Omega)
$$
is continuous for all $0\le s<s_0$.  
\end{thm}

\noindent
{\bf Remark.}
The sharp value of $s_0$ is not known, and most likely the exact
determination of such a value might prove a very difficult task.
The two sources \cite{BeCh} and \cite{Kohn-exp} present completely
different approaches and descriptions of $s_0$, that is of the range
$[0,s_0)$ for which $P$ is bounded on $H^s$, with $s\in[0,s_0)$.
In \cite{BeCh} it is proved that such a range is at least 
$[0,\delta/2)$, i.e. they show that $s_0\ge \delta/2$,  while in 
\cite{Kohn-exp} the parameter $s_0$ is not so explicit, but it 
tends to infinity 
as $\delta\rightarrow1$.
The value found in \cite{BeCh} has the advantage of providing an
explicit lower bound for the regularity of the Bergman projection on a
given domain, while the value given in \cite{Kohn-exp} is sharp in the
sense given by Boas and Straube's result Theorem \ref{psh-cond-R}.
\medskip

The domain $\sW$ serves as an example that the exponent $\d$ may
be arbitrarily small.  To illustrate this point, the following result
is essentially proved in 
\cite{DFo1}. Here we add the precise estimate that such an exponent $\d$
is less than the value $\nu$.
\begin{thm}\label{DF-exponent}	   \sl
Let $\d_0>0$ be fixed.  Then there exists $\mu_0>0$ such that
for all $\mu\ge\mu_0$ the following holds.  If $\widetilde\varrho$
is a defining function for 
$\sW=\sW_\mu$, 
with $\mu\ge\mu_0$ and $\d>0$ is such that
$-(-\widetilde\varrho)^\d$ is a bounded plurisubharmonic exhaustion
function for $\sW$, then $\d<\d_0$.

More precisely, we show that, in the notation above, $\d<\nu=\pi/2\mu$.
\end{thm}
\proof
We may assume that $\widetilde\rho=h\rho$, where $\rho=\rho_\mu$ is
defined 
in (\ref{def-rho}) and $h$ is a smooth positive function on
$\overline{\sW}$. Then, by hypothesis $-h^\d(-\rho)^\d$ is
strictly plurisubharmonic on $\sW$.

Let 
\begin{align*}
\sigma(z_1,z_2)
& = -\frac{1}{2\pi} \int_0^{2\pi} h^\d(z_1,e^{i\theta}z_2)
\bigl(-\rho(z_1,e^{i\theta}z_2)\bigr)^\d \, d\theta \\
& = -\frac{1}{2\pi} \int_0^{2\pi} h^\d(z_1,e^{i\theta}z_2)\, d\theta 
\bigl(-\rho(z_1,z_2)\bigr)^\d \\
& = -\widetilde h(z_1,z_2) \bigl(-\rho(z_1,z_2)\bigr)^\d \, .
\end{align*}
 Obviously, $\sigma$ is also
strictly plurisubharmonic on $\sW$, and $\widetilde h$ is strictly
positive and smooth on $\sW$.  We can also write
$\widetilde h(z_1,z_2) = h^{\#}(z_1,|z_2|^2) $, where $h^{\#}$ is
defined
for $(z_1,t)\in\bC\times\bR^+$ such that if $|z_2|^2=t$ then
$(z_1,z_2)\in\sW$. For simplicity of notation, we rename such a
function $h$ again.

Thus we have that
$$
\sigma(z_1,z_2)=
-h(z_1,|z_2|^2) \bigl(-\rho(z_1,z_2)\bigr)^\d 
$$
is strictly plurisubharmonic on $\sW$.  

Now consider the points in $\sW$ of the form $p = (z_1,z_2)= (\ve
e^{i\log|z_2|^2},z_2)$ with $e^{-\mu/2}\le |z_2|\le e^{\mu/2}$. 
For these points one has that
$$
\p\rho(p) =\bigl( (1-\varepsilon)e^{i\log|z_2|^2},0\bigr)\ .
$$
A straightforward computation shows that, at such points $p
\equiv (\ve e^{i\log|z_2|^2},z_2)$ the Levi
form  $\mathcal L_\sigma$ of $\sigma$ 
calculated at vectors $v=(v_1,v_2)\in\bC^2$ equals (all the functions
are evaluated at the points $p$
and we write $\zeta$ in place of $ e^{i\log|z_2|^2}$)
\begin{align*}
\mathcal L_\sigma \bigl(p;(v_1,v_2)\bigr)
& = \ve^{\d-2} (2-\ve)^{\d-2} \Bigl\{ 
(2-\ve)\Bigl( -\ve^2 (2-\ve) \p_{z_1\bar z_1}^2 h 
+2\delta\ve(1-\ve)\Re \bigl(\zeta \p_{z_1}h\bigr) +\delta\ve h\\
& \qquad\qquad\qquad\qquad
+\delta(1-\delta)\frac{(1-\ve)^2}{2-\ve}h \Bigr) |v_1|^2 \\
& \qquad +2\ve(2-\ve) \Re \Bigl[\bigl( -\ve(2-\ve)\p_{z_1\bar z_2}^2 h 
+\delta(1-\ve)\p_{z_2}h +\delta\frac{i\zeta}{z_2} h\bigr) v_1\bar v_2
\Bigr] \\ 
& \qquad 
+ \d^2(2-\d) \Bigl( -(2-\d)\p_{z_2\bar z_2}^2 h
+\frac{2\delta}{|z_2|^2}h \Bigr) |v_2|^2  \Bigr\} \ .
\end{align*}
Next, we evaluate the above Levi form at vectors of the form
$(v_1,v_2)=(u_1, \ve u_2)$.  Making the obvious simplification, we see
that the necessary condition in order for $\sigma$ to be strictly
plurisubharmonic is that 
and $0<\ve<1$ and, for all $(u_1,u_2)\in\bC^2$,
\begin{multline*}
\Bigl( -\ve^2 (2-\ve) \p_{z_1\bar z_1}^2 h 
+2\delta\ve(1-\ve)\Re \bigl(\zeta \p_{z_2}h\bigr) +\delta\ve h
+\delta(1-\delta)\frac{(1-\ve)^2}{2-\ve}h \Bigr) |u_1|^2 
\\
+ 2\Re \Bigl[\bigl( -\ve(2-\ve)\p_{z_1\bar z_2}^2 h 
+\delta(1-\ve)\p_{z_2}h +\delta\frac{i\zeta}{z_2} h\bigr) u_1\bar u_2
\Bigr] 
+  \Bigl( -(2-\ve)\p_{z_2\bar z_2}^2 h
+\frac{2\delta}{|z_2|^2}h \Bigr) |u_2|^2 \ge0\ .
\end{multline*}

Since $h\in C^\infty(\overline{\sW})$ this inequality must hold also
for $\ve=0$ and $(0,z_2)\in\mathcal A$.  Then we have
\begin{equation}
\Bigl( \frac12 \delta(1-\delta) h \Bigr) |u_1|^2 
+ 2\Re \Bigl[\bigl( \delta \zeta \p_{z_2}h 
+\delta\frac{i\zeta}{z_2} h\bigr) u_1\bar u_2
\Bigr] 
+  \Bigl( -(2-\ve)\p_{z_2\bar z_2}^2 h
+\frac{2\delta}{|z_2|^2} \Bigr) |u_2|^2 \ge0\ . 
\label{not-understood}
\end{equation} 

Next, we substitute for $h$ the function $\widetilde h$ defined on
$\bC\times\bR^+$  such that $h(z_1,z_2) = 
\widetilde h(z_1,|z_2|^2)$.  Then 
$$
\p_{z_2} h(0,z_2) = \bar z_2 \p_t
\widetilde h (0,|z_2|^2)\qquad\text{and}\qquad
\p_{z_2\bar z_2}^2 h(0,z_2) = |z_2|^2 \p_t^2 \widetilde h(0,|z_2|^2)
+  \p_t \widetilde h (0,|z_2|^2)\ .
$$
Plugging these into (\ref{not-understood}) we then obtain the
differential inequality for the function $\widetilde h$:
$$
 \frac12 \delta(1-\delta) \widetilde h  |u_1|^2 
+ 2\Re \Bigl[\delta\zeta\bigl( \p_{t}\widetilde h 
+\frac{i}{|z_2|^2} \widetilde h\bigr) u_1\bar u_2
\Bigr] 
+  \Bigl( -2|z_2|^2 \p_{t}^2 \widetilde  h -2\p_t \widetilde h 
+\frac{2\delta}{|z_2|^2} \widetilde h \Bigr) |u_2|^2 \ge0\ 
$$
for all $(u_1,u_2)\in\mathbf C^2$, $e^{-\mu/2}\le|z_2|\le e^{\mu/2}$
(and the function  $\widetilde h$ being evaluated at the points
$(0,|z_2|^2)$). Now if we choose $(u_1,u_2)$ of the form
$(2e^{i\theta}/|z_2|,1)$ in such a way that the second term in the
above display becomes non-positive, we obtain that the function
$\sigma$ is plurisubharmonic only if 
$$
 \frac{\delta(1-\delta)}{|z_2|^2}  \widetilde h  
-2\delta\bigl( (\p_{t}\widetilde h)^2
+\frac{\widetilde h^2}{|z_2|^4} 
\bigr)^{1/2} 
-|z_2|^2 \p_{t}^2 \widetilde h -\p_t \widetilde h
+\frac{\delta}{|z_2|^2} \widetilde h  \ge0
$$
which in turns gives
\begin{equation}\label{simplified-differential-inequality}
-\delta^2  \widetilde h  
-|z_2|^4 \p_{t}^2 \widetilde h -|z_2|^2 \p_t \widetilde h 
\ge0
\end{equation}
for all points $(0,|z_2|^2)$ with $e^{-\mu/2}\le|z_2|\le e^{\mu/2}$.

We now set $g(s)=\widetilde h(0,e^s)$ for $s\in[-\mu,\mu]$.  Notice
that $|z_2|^2=e^s$ and that
$$
g'=e^s \p_t \widetilde h\qquad\text{and}\qquad
g''=  e^s \p_t \widetilde h +e^{2s} \p_t^2 \widetilde h\ .
$$
From (\ref{simplified-differential-inequality}) we obtain the 
differential inequality
$$
g''+\delta^2 g
\le0\ ,
$$
for $s\in[-\mu,\mu]$, where $g$ is a smooth strictly positive
function. 
From the strict positivity of $g$ it follows that, for all
$0<\delta'<\delta$, it must be that
$$
g''+{\delta'}^2 g <0\ ,
$$
again for all $s\in[-\mu,\mu]$.  Setting $\widetilde g(s)=g(s/\delta')$
the differential inequality above can be re-written as
$$
\widetilde g'' +\widetilde g<0
$$
for all $s\in[-\mu\delta',\mu\delta']$. Finally, by  translation
(calling the new function $g$ again), i.e. setting
$g(s)=\widetilde g(s+\mu\delta')$, 
 we obtain that
\begin{equation}\label{final-differential-inequality}
g''+g<0
\end{equation}
for a smooth strictly positive function $g$, 
for all $s\in[0,2\mu\delta']$.  

We now claim that there exists a smooth strictly positive function
$\varphi$ such that 
\begin{equation}\label{double-final-differential-inequality}
\varphi''+\varphi<0\qquad\text{ and}\qquad \varphi'<0 
\end{equation}
$s\in[0,\mu\delta']$.  For notice that if $g$ as above is such that
$g'(a)<0$, then $g'(s)<0$ for $s\in[a,2\mu\delta']$, while,  if instead
$g'(a)\ge0$, then $g'(s)>0$ for $s\in[0,a)$, since $g''<0$ on
$[0,2\mu\delta']$. In this latter case, making the substitution
$s\mapsto 2\mu\delta'-s$ that preserves
(\ref{final-differential-inequality}), we obtain a function with
negative derivative on $[a,2\mu\delta')$.  By the arbitrariness of
$\delta'<\delta$ we establish the claim.

Now, the argument at the end of the proof of Theorem 6 in \cite{DFo1}
shows that the differential inequalities
(\ref{double-final-differential-inequality}) above are possible only if
$\mu\delta' < \pi/2$, i.e. 
$$
\delta'<\frac{\pi}{2\mu}=\nu\ .
$$
This proves the result.
\endpf

\bigskip
\section{Non-Smooth Versions of the Worm Domain}\label{3}

In order to perform certain analyses on $\sW$ some
simplifications  of the domain turn out to be particularly useful.

In the first instance, one can simplify the expression of the 
defining function $\rho$ for $\sW$ by taking $\eta$ to be
the characteristic function of the interval $[-\mu,\mu]$. 
This has the effect of truncating the two caps and
destroying in part the smoothness of the boundary.  
Precisely, we can
define
\begin{equation}\label{Omega-beta-prime}
\sW' =
\Big\{(z_1, z_2)\in\bC^2:\,  
\big|z_1 - e^{i\log |z_2|^2} \big|^2 < 1,\  
\big|\log|z_2|^2\big|<\mu \Big\}
\, .
\end{equation}

We remark that
$\sW'$ is a bounded, pseudoconvex domain with boundary that
is $C^\infty$ except at points that satisfy
\begin{itemize}
\item[{\bf (i)}] $|z_2| = e^{\mu/2}$ and $|z_1 - e^{-i\log|z_2|^2}| = 1$;
\item[{\bf (ii)}]  $|z_2| = e^{-\mu/2}$ 
and $|z - e^{-i\log|z_2|^2}| = 1$.\medskip
\end{itemize}

Of interest are also two non-smooth, {\it unbounded} worms.
Here, in order to be consistent with the results obtained in
\cite{KrPe}, 
we change the notation a bit.  (In practice, we set
$\mu=\beta-\pi/2$.)

For $\b>\pi/2$ we define
\begin{equation}\label{D-beta}
D_\b = \Big\{\z \in \bC^2: 
\Re \bigl(\z_1 e^{-i\log |\z_2|^2}
\bigr)>0 ,\, \big|\log |\z_2|^2\big| < 
\b - \frac{\pi}{2} \Big\} 
\end{equation}
and
\begin{equation}\label{D-beta-prime}
D'_\b = \Big\{z 
\in \bC^2: \big|\Im z_1 - \log |z_2|^2\big| 
< \frac{\pi}{2},\,  |\log |z_2|^2| <
\b - \frac{\pi}{2} \Big\} \, . 
\end{equation}

It should be noted that these  latter two domains 
are biholomorphically equivalent via the mapping
\begin{equation}\label{eq:biholo}
(z_1,z_2)\ni D_\b' \mapsto (e^{z_1},z_2)\ni D_\b\, .
\end{equation}
Neither of these domains is bounded.
Moreover, these domains are {\it not} smoothly
bounded.   Each boundary is only Lipschitz, and, in particular, their
boundaries are Levi flat. 

We notice in passing that the slices of $D_\b$, for each fixed
$\zeta_2$, are halfplanes in the variable $\zeta_1$.  Likewise
the slices of $D'_\b$, for each fixed $\zeta_2$, are strips in
the variable $\zeta_1$.

The geometries of these domains are rather different from that of the
smooth worm 
$\sW$, which has smooth boundary and all boundary points, except those
on a singular 
annulus $(0, e^{i\log |z_2|^2})$ in the boundary, are strongly pseudoconvex.
However our worm domain $D_\b$ is actually a model for the smoothly
bounded $\sW$ (see, for instance, \cite{Ba3}), and it can be expected
that phenomena that are true on $D_\b$ or $D'_\beta$ will in fact hold
on $\sW$ as 
well. We will say more about this symbiotic relationship below.

We now illustrate a first application of these non-smooth domains in
the analysis of $\sW$.  We begin with the main result of \cite{Ki}.

\begin{thm}\label{Kiselman-thm}  \sl
Let $\sW'$ be as above.
Then, there is a function $f \in C^\infty(\overline{\sW'})$ 
such
that 
its Bergman projection $Pf$ is not H\"{o}lder continuous of any
positive order on 
$\overline{\sW'}$.
\end{thm}

Following Kiselman \cite{Ki}, we now describe an outline of the proof
of this theorem.  The steps are as follows:
\begin{enumerate}
\item[{\bf (a)}]  
We construct a subspace $C^+(\sW')$ of $L^2(\sW')$ which contains 
all the H\"{o}lder continuous functions on $\overline{\sW'}$.
\item[{\bf (b)}]  We construct a linear functional $T$ whose values are
obtained as holomorphic extensions of inner products $\langle f,
g_\alpha \rangle$ for certain elements $g_\alpha$ of the Bergman
space.   
That is to say, for a fixed $f \in C^+(\sW')$, we define a
holomorphic function of 
the complex variable $\alpha$ by $\Phi(\alpha) = \langle f, g_\alpha
\rangle$; here 
$\Re \alpha > -1$.  We set $F(f) = \Phi(-2)$. 
\item[{\bf (c)}]  We show that if $f$ and $Pf$ both belong to
  $C^+(\sW')$ then $T(Pf) = T(f)$; in 
particular, $f - Pf$ is orthogonal to ${\mathcal O}^2(\sW')$, hence
$f - Pf$ is 
orthogonal to the $g_\alpha$.
\item[{\bf (d)}]  We show that $T(f) = 0$ if $f$ is in $C^+(\sW')$
  and holomorphic. 
\item[{\bf (e)}]  We show that $T$ is not identically zero on
  $C^+(\sW')$.  Specifically, 
$C^+(\sW')$ contains $C^\infty(\overline{\sW'})$ and $T$ is not zero on
$C^\infty(\overline{\sW'})$.
\item[{\bf (f)}]  We finish the proof by taking an $f \in
  C^\infty(\overline{\sW'})$ with 
$T(f) \ne 0$.  If $Pf$ belongs to $C^+(\sW')$ then steps {\bf (a)}
and {\bf (c)} tells us that 
$T(Pf) = T(f) \ne 0$.  That contradicts {\bf (b)}.
\end{enumerate}

Kiselman's work was pioneering in that it put the worm domain at the
forefront for examples that bear on Condition $R$ and the regularity
of the $\dbar$ problem. 

The ``worm'' that 
Kiselman studies does
{\it not} have smooth boundary.  
Yum-Tong Siu 
\cite{Siu} 
later proved a version of
Kiselman's theorem 
on the smooth worm.  His result is:

\begin{thm}\label{siu-thm}   \sl
For a suitable version of the smooth worm $\sW$, there is
a function $f \in C^\infty(\overline{\sW})$ such that
the Bergman projection $Pf$ is not H\"{o}lder continuous 
of any positive order on $\sW$.
\end{thm}

Siu's proof is quite intricate, and involves an argument with de Rham
cohomology 
to show that caps may be added to Kiselman's domain to make it into a
smooth worm. 
\bigskip

\section{Irregularity of the Bergman Projection}\label{4}

We begin this section by discussing  the proof
of Barrett's result  
Theorem \ref{Barrett-thm} \cite{Ba1}.  Now let us describe these ideas
in some detail.  We begin with some of Kiselman's main ideas.

Let the Bergman space ${\cal H}=A^2$ be the collection of
holomorphic functions that are square integrable with
respect to Lebesgue volume measure $dV$ on a fixed domain.
Following Kiselman \cite{Ki} and Barrett \cite{Ba2}, using the
rotational invariance in the $z_2$-variable,  
we decompose the Bergman space for the domains $D_\b$ and $D_\b'$  as
follows. 
Using the rotational invariance in $z_2$ and 
elementary Fourier series, each $f \in{\cal H}$  can be  written as 
$$
f = \sum_{j=-\infty}^\infty f_j \, ,
$$
where each $f_j$ is holomorphic and satisfies 
$f_j(z_1, e^{i\theta} z_2) = e^{ij\theta} f(z_1, z_2)$ 
for $\theta$ real.
In fact such an $f_j$ must have the form
$$
f_j(z_1, z_2) = g_j(z_1, |z_2|) z_2^j \, ,
$$
where $g_j$ is holomorphic in $z_1$ and locally constant in $z_2$.

Therefore
$$
{\cal H} = \bigoplus_{j\in\ZZ} {\cal H}^j\, ,
$$
where
$$
{\cal H}^j = \bigl \{f \in L^2: f \ \hbox{is holomorphic and} \ 
          f(w_1, e^{i\theta} w_2) = e^{ij\theta} f(w_1, w_2) \bigr \} \, .
$$

If $K$ is the Bergman kernel for ${\cal H}$ and $K_j$ the
Bergman kernel for ${\cal H}^j$, then we may write
$$
K = \sum_{j=-\infty}^\infty K_j \, .
$$

Notice that, by the invariance property of ${\cal H}^j$, with $z =
(z_1, z_2)$ and $w = (w_1, w_2)$, we have that 
$$
K_j(z,w) = H_j(z_1,
w_1) z_2^j \overline{w}_2^j \,.  
$$
Our job, then, is to calculate each
$H_j$, and thereby each $K_j$.  The first step of this calculation is 
already done in \cite{Ba2}.  We outline the calculation here for the
sake of completeness. 
\begin{proposition}\label{H-j}  \sl
Let $\b>\pi/2$.  Then 
\begin{equation}\label{H-j-equation}
H_j(z_1,w_1) = \frac{1}{2\pi} \int_{-\infty}^\infty
   \frac{ e^{i(z_1 - \bar w_1)\xi} 
       \xi \left (\xi - \frac{j+1}{2} \right ) }{\sinh(\pi \xi)
        \sinh \left [(2\b - \pi) \left (\xi - \frac{j+1}{2}
         \right ) \right ] } 
        \, d\xi \, .
\end{equation}
\end{proposition}

The papers \cite{Ki} and [Ba2] calculate and analyze only the Bergman
kernel for ${\cal H}^{-1}$ (i.e., the Hilbert subspace with index $j = -1$).  
This is attractive to do because
certain ``resonances'' cause cancellations that make the
calculations tractable when $j = -1$.  One of the main thrusts
of the work \cite{KrPe} is to perform the more difficult
calculations for all $j$ and then to sum them over $j$. 
\medskip

\proof
We begin by following the calculations in \cite{Ki} and \cite{Ba2} in
order to get our hands on the Bergman kernels of the ${\cal H}^j$.
Let $f_j \in {\cal H}^j$ and fix $w_2$.  Then
$f_j(w_1, w_2) = h_j(w_1)  w_2^j$ (where we of
course take into account the local independence of $h_j$ from $w_2$).
Now, writing $w_1 = x + iy$, $w_2 = re^{i\theta}$, and
then making the change of variables  $\log r^2 = s$,
we have
\begin{eqnarray}
\|f_j\|_{{\cal H}}^2 
& = & \int_{D'_\b} |h_j(w_1)|^2 |w_2|^{2j} \, dV(w) \notag\\
& = & \int_{-\infty}^\infty \int_{|y - \log r^2| < \frac{\pi}{2}}
  2\pi |h_j(x+iy)|^2  \int_{|\log r^2| < \b - \frac{\pi}{2}}
    r^{2j+1} \, dr dy dx \notag\\
& = & \pi \int_\bR \int_{|y - s| < \frac{\pi}{2}} |h_j(x + iy)|^2 
    \int_{|s| < \b - \frac{\pi}{2}} e^{s(j+1)} \, ds dy dx \notag\\
& = & \pi \int_{|y| < \b,\,  x \in \bR} |h_j(x + iy)|^2 
   \int_{-\infty}^\infty e^{(j+1)s} \chi_{\pi/2} (y - s) \chi_{\b -
     \pi/2} (s) \,  
     ds dx dy \notag\\
& = & \int_{S_\b} |h_j(w_1)|^2 \biggl ( \chi_{\pi/2} * \bigl [
   e^{(j+1)( \, \cdot \, )} \chi_{\b - \pi/2}( \, \cdot \, ) \bigr ]
   \biggr ) (y) \, 
   dx dy \, ;  \label{star}
\end{eqnarray}
here we have set 
$$
S_\b= \{x+iy \in \bC:  |y | < \b\}
$$
and used the notation  
$$
\chi_\alpha(y) = \left \{ \begin{array}{lcl}
                         1 & \hbox{if} & |y| < \alpha \, , \\
                         0 & \hbox{if} & |y| \geq \alpha \, .
                       \end{array}
                 \right.  
$$

For $\b > \frac{\pi}{2}$, we now set
$$
\lambda_j(y) = \biggl ( \chi_{\pi/2} * \bigl [ 
   e^{(j+1)( \, \cdot \, )} \chi_{\b - \pi/2} ( \, \cdot \, )\bigr ]
   \biggr ) (y) \, 
   . 
$$
So line (\ref{star}) equals 
$$
 \int_{S_\b} |h_j(w_1)|^2 \lambda_j(y) \, dx dy \, .
$$
Thus we have shown that, if $f_j \in {\cal H}^j$,
$f_j = h_j(w_1) w_2^j$, then
$$
\|f_j\|_{{\cal H}}^2 
= \int_{S_\b} |h_j(w_1)|^2 \lambda_j(y) \, dx dy \, .
$$

Now let $\varphi \in A^2(S_\b, \lambda_j \, dA)$.  That is,
$\varphi$ is square-integrable on $S_\b$ with respect to the measure
$\lambda_j \, dA$ (here $dA = dx dy$ is two-dimensional area
measure). 
Note that $\lambda_j$ depends only on the single variable $y$.  Let
$\widetilde{\varphi}$ denote the partial Fourier transform 
of $\varphi(x + iy)$ in the $x$-variable.
Then (by standard Littlewood-Paley theory)
$$
\widetilde{\varphi}(\xi, y) 
=  \int \varphi(x + iy) e^{-ix\xi} \, dx  
=  e^{-y\xi} \widetilde{\varphi}_0 (\xi) \, , \\
$$
where $\varphi_0(x) = \varphi(x + i0)$.  Therefore,
denoting by $B_\b = B_\b^{(j)}$  the Bergman kernel for the strip $S_\b$
with respect to the weight $\lambda_j$ and writing
$\omega=s+it$ and denoting by $\xi$ the variable dual to $s$, we have   
\begin{eqnarray*}
\int_\bR \widetilde{\varphi}_0 (\xi) e^{i\z \xi} \, d\xi 
& = & 2\pi \varphi(\z)  =  2\pi \int_{S_\b} \varphi(\omega) B_\b(\z, \omega)
       \lambda_j(\Im \omega) \, dA(\omega)  \\
& = & \int_{-\b}^\b \int_\bR \widetilde{\varphi} (\xi, t)
       \widetilde{B}_\b (\z, (\xi, t)) \lambda_j(t) \,
           d\xi dt \\
& = & \int_\bR \widetilde{B}_\b (\z, (\xi, 0)) 
\int_{-\b}^\b \widetilde{\varphi}_0 (\xi)
         e^{-2\xi  t} 
               \lambda_j(t) \, dt \, d\xi \, .
\end{eqnarray*}
Notice that there is a factor of $e^{-\xi t}$ from each
of the Fourier transform functions in the integrand.

This gives a formula
for $\widetilde{B}_\b$:
$$
\widetilde{B}_\b (\z, (\xi, 0)) = \frac{e^{i\z\xi}}{\int_{-\b}^\b
       e^{-2t\xi} \lambda_j(t) \, dt} = 
       \frac{e^{i\z\xi}}{\widehat{\lambda}_j(-2i\xi)} \, .
$$

Amalgamating all our notation, and using the fact that
the (Hermitian) diagonal in $\bC^2$ is a set of
determinacy, we find that
$$
B_\b(z,w) = \frac{1}{2\pi} \int_\bR \frac{e^{i(z -
    \overline{w})\xi}}{\widehat{\lambda}_j(-2i\xi)} \, d\xi \, . 
$$
But of course
$(\chi_{\pi/2})\, {}\widehat{\mathstrut} \, 
(\xi)  = (e^{i\xi \pi/2} - e^{-i\xi \pi/2})/\xi$,
so that
$$
(\chi_{\pi/2}) \, \, {}\widehat{\mathstrut} \, \, (-2\xi i) 
 = 
 \frac{1}{\xi} \sinh (\pi \xi ) \, .
$$
Furthermore,
$$
\left ( e^{(j+1)s} \chi_{\b - \frac{\pi}{2}} (s) \right ) 
 \, {}\widehat{\mathstrut} \, =  \frac{\sinh \bigl( (2\b - \pi)
            \bigl(\xi - \frac{j+1}{2} \bigl) 
            \bigr)}{\xi - \frac{j+1}{2}}  \, .
$$
Thus 
$$
\widehat{\lambda}_j(-2i\xi) = \frac{\sinh ( \pi \xi)
     \sinh((2\b - \pi)(\xi - \frac{j+1}{2}))}{\xi(\xi 
         - \frac{j+1}{2} )} \, 
$$
and
$$
\frac{1}{\widehat{\lambda}_j(-2i\xi)} = \frac{\xi\left (
     \xi -\frac{j+1}{2} \right )  }{\sinh ( \pi \xi)
     \sinh((2\b - \pi)(\xi - \frac{j+1}{2}))} \, .
$$
In conclusion, 
$$
H_j(z_1, w_1) = \frac{1}{2\pi} \int_{-\infty}^\infty
   \frac{ \left (e^{i(w_1 - \overline{z}_1)\xi} \right )
       \xi  \left (\xi - \frac{j+1}{2} \right ) }{\sinh(\pi \xi)
       \sinh \left ((2\b - \pi) \left (\xi - \frac{j+1}{2}
         \right ) \right ) } 
        \, d\xi \, ,
$$
thus proving (\ref{H-j-equation}). 
\endpf

At this point we sketch the proof of the main result of Barrett in
\cite{Ba3}. 

\proof[Sketch of the proof of Theorem \ref{Barrett-thm}]
The proof starts from the observation that the Bergman projection
$\mathcal P$ on $\mathcal W$ preserves each $\mathcal H^j$.  Therefore
in order to show that $\mathcal P$ is not continuous on $H^s$, for
some $s$, it suffices to show that 
$\mathcal P$ fails to be  continuous when restricted to some $\mathcal
H^j$.

The first step is to calculate the asymptotic expression for the
kernel when $j=-1$.  Recall that we are working on the non-smooth
domain $D_\beta'$.
Using the method of contour integral it is not difficult
to obtain that 
$$
K'_{-1}(z,w) = \bigl(
e^{-\nu_\beta|z_1-\bar w_1|} +\mathcal O (e^{-\nu|\Re z_1-\Re w_1|})
\bigr) \cdot(z_2\bar w_2)^{-1}
$$
as $|\Re z_1-\Re w_1|\raw+\infty$, 
uniformly in all closed strips $\{|\Im z_1|,\, |\Im w_1| \le ...\}$,
with $\nu > \nu_b$.  

By applying the biholomorphic transformation (\ref{eq:biholo}) one
obtains an asymptotic expression for the kernel $K_{-1}$ relative to
the domain $D_\beta$:
$$
K_{-1}(\z,\om) = (|\z_1||\om_1|)^{-1} \cdot\Bigl(
\frac{|\om_1|^{\nu_\b}}{|\z_1|^{\nu_\b}}
+\mathcal O( |\om_1|^{\nu_\b}/|\z_1|^{\nu_\b})^{-\nu} \Bigr)  
\cdot (\z_2\bar\om_2)^{-1}\ ,
$$
with $\nu>\nu_\b$, 
as $|\z_1|-|\om_1|\raw0^+$, 
The proof of these two assertions can be found in \cite{Ba3} (or see
\cite{CheS}).

The next step is a direct calculation to show that
$K_{-1}(\cdot,w)\not\in H^s(D_\beta)$ for $s\ge \nu_\b$.   This
assertion is proved by using the characterization of Sobolev norms for
holomorphic functions on a domain $\Omega$:  For $-1/2<t<1/2$, $m$ a
non-negative integer, the norm
$$
\sum_{|\alpha|\le m} \big\| |\rho|^t \p^\alpha_{z} h \big\|_{L^2(\Omega)}
$$
is equivalent to the $H^{m-t}$-norm of the holomorphic function $h$.
The proof of such a characterization can be found in \cite{Lig2}.

Next, one notice that the reproducing kernel $K_{-1}(\cdot,w)$ can be
written as the projection of a radially symmetric smooth cut-off
function $\chi$, translated at $w$.  That is, if we denote by $P_{-1}$
the projection relative to the subspace $\mathcal H^{-1}$, then
$$
K_{-1}(\cdot,w) = P_{-1} \bigl(\chi(\cdot-w)\bigr)\, .
$$
Therefore, since $K_{-1}(\cdot,w)\not\in H^s(D_\beta)$ for $s\ge
\nu_\b$, then $P_{-1}$, and therefore $P_{D_\b}$ is not continuous on
$H^s(D_\b)$.  \medskip

The final step of the proof is to transfer this negative result from
$D_\b$ to $\mathcal W$.  This is achieved by an exhaustion
argument. We adapt this kind of argument to obtain a negative result
in the $L^p$-norm in the proof of Theorem \ref{straube-thm} and we do
not repeat the argument here. 
\endpf				       

\bigskip

\section{Failure of Global Hypoellipticity and Condition
  $R$}\label{ConditionR}

In order to discuss the failure of Condition $R$ on the
Diederich-Forn\ae ss worm domain, we recall the basic facts about the
$\dbar$-Neumann problem. 

Let $\Omega \subset \bC^n$ be a bounded domain with smooth
boundary and let $\rho$ be a smooth defining function for $\Omega$.
The $\dbar$-Neumann problem on $\Omega$ is a boundary value problem
for the elliptic partial differential operator
$$
\BoxOpTwo = \dbar
\dbar^* + \dbar^* \dbar\ .
$$
Here $\dbar^*$ denotes the $L^2$-Hilbert space  adjoint of the
(unbounded) operators $\dbar$.
In order to apply $\BoxOpTwo$ to a form or current $u$ one
needs to require that $u,\dbar u\in\dom(\dbar^*)$.
These conditions translate into two differential equations on the
boundary for $u$
 and the are called the {\it two $\dbar$-Neumann boundary conditions} 
(see \cite{FoKo} or \cite{Treves}). 
These equations are
\begin{equation}\label{Neumann-bndry-conditions}
u\lrcorner\dbar\rho =0\, ,\qquad\text{and}\quad
\dbar u\lrcorner\dbar\rho =0\, ,\qquad\text{on}\ \p\Omega\ .
\end{equation}
Thus the equation $\BoxOpTwo u=f$ becomes a boundary
value problem. 
\begin{equation}\label{BVP}
\left\{
\begin{array}{lcl}
\BoxOpTwo u & = & f \quad \hbox{on} \ \Omega   \\
u\lrcorner\dbar\rho \, ,
\dbar u\lrcorner\dbar\rho & = & 0 \quad\text{on}\ \p\Omega\ .
\end{array}
\right.
\end{equation}
This is an equation defined on {\it forms}.  
The significant problem is for $(0,1)$-forms, and we restrict to this
case in the present discussion.

It follows from H\"ormanders' original paper on the solution of the
$\dbar$-equation \cite{Hormander-Acta-paper}  that the $\dbar$-Neumann
problem is always solvable on a smoothly bounded pseudoconvex domain
$\Omega$ in
$\bC^n$ for any data $f\in L^2(\Omega)$. We denote by $N$---the
Neumann operator---such a  solution operator.  Moreover, $N$ turns
out to be continuous in the $L^2$-topology:
$$
\|Nu\|_{L^2} \le c\|u\|_{L^2}\ .
$$
An important
formula of Kohn says 
that 
$$
P = I - \dbar^* N \dbar \, .
$$
The proof of this is a formal calculation---see \cite{Kr2}.  
Important work by Boas and Straube \cite{BoS2}
essentially established that
the Neumann operator $N$ has a certain regularity (that is, it 
maps some Sobolev space $H^s$ to itself, for instance) if and only if
$P$ will have the same regularity property.  In particular if $N$ is
continuous on a Sobolev space $H^s$ for 
some $s>0$ (of $(0,1)$-forms), then the
Bergman projection $P$ is continuous on the same Sobolev space $H^s$
(of functions).

Such regularity is well
known to hold on strongly pseudoconvex domains (\cite{FoKo}, \cite{Kr2}).
In addition, Catlin proved a similar regularity result on finite
type domains (see \cite{Cat1}, \cite{Cat2}, \cite{Kr1}).\medskip

Michael Christ's milestone result \cite{Chr1} has proved to be
of central importance for the field.  It demonstrates concretely
the seminal role of the worm, and points to future directions
for research.  Certainly the research program being described here,
including the calculations in \cite{KrPe}, is inspired by Christ's work.

Christ's work is primarily concerned with 
global regularity, or global hypoellipticity.  A partial differential
operator $L$ is said to be {\it globally hypoelliptic} if, whenever 
$L u = f$ and $f$ is globally $C^\infty$, then $u$ is globally $C^\infty$.
We measure regularity, here and in what follows, using the standard
Sobolev spaces $H^s$,  
$0<s<\infty$  see \cite{Kr2}, \cite{Hor1}).

Christ's proof of the failure of global hypoellipticity is a highly
complex and recondite calculation with pseudodifferential operators.
We cannot replicate it here.  But the ideas are so important that
we feel it worthwhile to outline his argument.  We owe a debt to
the elegant and informative paper \cite{Chr2} for these ideas.
\medskip

As a 
boundary value problem for an elliptic operator, the $\dbar$-Neumann
problem  may be treated by Cald\'{e}ron's method of reduction to
a pseudodifferential equation on $\p\Omega$.  The sources
\cite{Hor2} and \cite{Treves} 
give full explanations of the classical
ideas about this reduction.  
In the more modern reference \cite{CNS} 
Chang, Nagel and Stein elaborate
the specific application of these ideas to the 
$\dbar$-Neumann problem in $\bC^2$. 
(Thus, in the remaining part of this discussion, $\Omega$ will denote a
smoothly bounded pseudoconvex domain in $\bC^2$.)
The upshot is that one reduces  the solution of the equation 
$\BoxOpTwo u=f$ to the solution of an equation 
$\BoxOpTwo^+ v=g$ on the boundary.  Here $u$ and $f$ are
$(0,1)$-forms, while
$v$ and $g$ are sections of a certain complex line bundle on
$\p\Omega$. (The fact that this bundle is 1-dimensional is a
consequence of the inclusion $\Omega\subset\bC^2$.) 

To be more explicit, the solution $u$ of (\ref{BVP}) can
be written as $u=\mathcal G f+ \mathcal R v$, where $\mathcal G$ is the
{\em Green operator} and $\mathcal R$ is the {\em Poisson
  operator}\footnote{Thus $\mathcal G$ is the solution operator for
  the 
elliptic boundary value problem $\BoxOpTwo(\mathcal G f) = f$ on 
$\Omega$ and
$\mathcal Gf = 0$ on $\p\Omega$, while 
$\mathcal R$ is the solution operator for
  the 
elliptic boundary value problem $\BoxOpTwo(\mathcal R v) = 0$ 
on $\Omega$ and $\mathcal R v=v$ on $\p\Omega$.}
for the operator $\BoxOpTwo$ and $v$ is chosen in such a way to
satisfy the boundary conditions.  In fact,
\begin{align*}
& \BoxOpTwo \bigl(\mathcal G f+ \mathcal R v\bigr) = f+0=
f \qquad\text{on}\
\Omega \\ 
& \bigl(\mathcal G f+ \mathcal R v\bigr)\lrcorner \dbar\rho =
 v \lrcorner \dbar\rho = 0 \qquad\text{on}\
\p\Omega \\ 
& \dbar \bigl(\mathcal G f+ \mathcal R v\bigr)\lrcorner \dbar\rho =
\dbar\mathcal G f\lrcorner \dbar\rho + 
\dbar v \lrcorner \dbar\rho = 0 
\qquad\text{on}\
\p \Omega \ .
\end{align*}
The section $v$ has two components, but one of these vanishes because
of the first $\dbar$-Neumann boundary condition.  
The second $\dbar$-Neumann
boundary condition may be written as an equation $\BoxOpTwo^+ v = g$
on $\p\Omega$,  
where $\BoxOpTwo^+$ is a pseudodifferential operator of order 1.  Also we
note that $g = (\dbar \mathcal G f \lrcorner \dbar \rho)$ restricted 
to $\p\Omega$.
\medskip

Christ's argument begins with a real-variable model for the 
$\dbar$-Neumann
problem that meshes well with the geometry of the boundary of the
worm domain $\mathcal W$.  

Let $M$ be the 2-torus $\mathbf T^2$ and
let $X, Y$ two 
smooth real vector 
fields on $M$.  Fix a coordinate patch $V_0$ in $M$ and suppose
that $V_0$ has been identified with 
$\{(x,t) \in (-2,2) \times (-2\d, -2\d)\} \subset \bR^2.$
Let $J = [-1,1] \times \{0\} \subset V_0$.

Call a piecewise smooth path $\gamma$ on $M$ {\it admissible}
if every tangent to $\gamma$ is in the span of $X, Y$.  Assume that 
\begin{enumerate}
\item[{\bf (i)}]  The vector fields $X, Y, [X,Y]$ span the tangent space to 
$M$ at every point of $M \setminus J$.
\item[{\bf (ii)}]  In $V_0$, $X \equiv \p_x$ and $Y \equiv
  b(x,t) \p_t$. 
\item[{\bf (iii)}]  For all $|x| \leq 1$ and $|t| \leq \d$, we have that
$b(x,t) = \alpha(x)t + {\cal O}(t^2)$, where $\alpha(x)$ is nowhere vanishing.
\end{enumerate}
It follows then that every pair $x, y \in M$ is connected by an
admissible path. 

\begin{thm}\label{5.2} \sl
With $X, Y, M$ as above, let $L$ be any partial differential operator
on $M$ of the form $L = - x^2 - Y^2 + a$, where $a \in C^\infty(M)$ and
\begin{equation}\label{eqno8}
\|u\|^2 \leq C \langle Lu, u \rangle  
\end{equation} 
for all $u \in C^2(M)$.  Then
$L$ is not globally regular in $C^\infty$.
\end{thm}

We note that our hypotheses, particularly inequality (\ref{eqno8}),
imply that $L$ has a well-defined inverse $L^{-1}$ which is a bounded
linear operator on $L^2(M)$.

The following theorem gives a more complete, and quantitative, version
of this result:

\begin{thm} \sl
Let $X, Y, M, L$ be as above.  Then $L$ has the following global
properties:
\begin{enumerate}
\item[{\bf (a)}]   There is a positive number $s_0$ such that, for
every $0 < s < s_0$, $L^{-1}$ preserves $H^s(M)$;
\item[{\bf (b)}]   For each $s > s_0$, $L^{-1}$ fails to map
$C^\infty(M)$ to $H^s(M)$;
\item[{\bf (c)}]  There is a sequence of values $s < r$ tending to
  infinity such that if $u \in H^s(M)$ satisfies $Lu \in H^r(M)$ then
  $u \in H^r$; 
\item[{\bf (d)}]  There are arbitrarily large values of $s$ with
a constant $C = C_s$ such that if $u \in H^s(M)$ is such that $Lu \in
H^s(M)$ then
\begin{equation}\label{eqno9}
\|u\|_{H^s} \leq C \|Lu\|_{H^s} \, .  
\end{equation}
\item[{\bf (e)}]
 For each value of $s$ as in part {\bf (d)}, $\{f \in H^s(M): L^{-1} f
 \in H^s(M)\}$ is a 
closed subspace of $H^s$ with finite codimension.
\end{enumerate}
\end{thm}

The proof of Theorem \ref{5.2} breaks into two parts.  The first part
consists of proving 
the {\it a priori} inequality (\ref{eqno8}).  The second part,
following ideas 
of Barrett in \cite{Ba3}, shows that, for any $s \geq s_0$, the
operator $L$  
cannot be exactly regular on $H^s(M)$.  We refer the reader to \cite{Chr1}
for the details.  Section 8 of \cite{Chr2} also provides a nice outline
of the analysis.

The next step is to reduce the analysis of the worm domain, as defined in
our Sections \ref{2} and \ref{3}, to the study of the manifold $M$ as
above.  With this 
idea in mind we 
set $\overline{L} = \dbar_b$ and $L$ its complex conjugate.
The characteristic variety\footnote{The characteristic variety of a
  pseudodifferential operator is the conic subset of the cotangent
  bundle on which its principal symbol vanishes.}
 of $\overline{L}$ is a real line bundle $\Sigma$ that splits smoothly
 as two rays: $\Sigma=\Sigma^+ \cup \Sigma^-$.

The principal symbol of $\BoxOpTwo^+$ vanishes only on 
$\Sigma^+$ that is half the characteristic variety 
We may compose
$\BoxOpTwo^+$ with an elliptic pseudodifferential operator of order
$+1$ to change 
$\BoxOpTwo^+$ to the form
\begin{equation}\label{15?}
{\cal L} = \overline{L}L + B_1 \overline{L} + B_2 L + B_3 
\end{equation}
microlocally in a conical neighborhood of $\Sigma^+$, where each
$B_j$ is a pseudodifferential operator with order not exceeding 0.
Since $\BoxOpTwo^+$ is elliptic on the complement of $\Sigma^+$, our
analysis may thus be microlocalized to a small conical neighborhood of
$\Sigma^+$. 

For a worm domain ${\cal W}$, there is circular symmetry in the second
variable.  This induces a natural action on functions and on forms
(as indicated in Section \ref{4}).   As indicated earlier, the
Hilbert space of square integrable $(0,k)$-forms has the
orthogonal decomposition $\oplus_j {\cal H}^j_k$.  The Bergman
projection and the Neumann operator preserve ${\cal H}^j_0$ and
${\cal H}^j_1$.  We now have the following key result:

\begin{proposition} \sl
Let ${\cal W}$ be the worm.  Then there is a discrete subset $S
\subset \bR^+$ 
such that, for each $s \not \in S$ and each $j \in \ZZ$, there is a
constant 
$C = C(s,j) < \infty$ such that, for each $(0,1)$ form $u \in {\cal
  H}^j_1 \cap C^\infty(\overline{\cal W})$ 
such that $N u \in C^\infty$, it holds that
$$
\|N u \|_{H^s({\cal W})} \leq C \cdot \|u \|_{H^s({\cal W})} \, .
$$
\end{proposition}

The operators ${\cal L}$, $\overline{L}$, $L$, $B_j$ in
(\ref{eqno9})  may be
constructed 
so as to commute with the circle action in the second variable, hence they
will preserve each ${\cal H}^j$.  In summary, for each $j$, the
action of ${\cal L}$ on ${\cal H}^j(\p {\cal W})$ may be identified
with the action of an operator ${\cal L}_j$ on $L^2(\p {\cal W}/S^1)$.

Of course $\p {\cal W}$ is 3-dimensional, hence $\p {\cal  W}/S^1$ 
is a real 2-dimensional manifold.  It is convenient to take coordinates
$(x, \theta, t)$ on $\p {\cal W}$ so that 
$$
z_2 = \exp(x + i\theta) \quad \hbox{and} \quad z_1 
= \exp(i2x)(e^{it} - 1) \, ;
$$
here $|\log |z_2|^2| \leq r$ and ${\cal L}_j$ takes the form
$\overline{L}L + B_1 \overline{L} + 
B_2 L + B_3$ (just as in (\ref{15?})!).  In this last formula,
$\overline{L}$ is a complex 
vector field which has the form $\overline{L} = \p_x + it
\alpha(t) \p_t$, where 
$|x| \leq r/2$, $\alpha(0) \ne 0$, and each $B_j$ is a classical
pseudodifferential operator 
of order not exceeding 0---depending on $j$ in a non-uniform manner.

We set $J = \{(x,t): |x| \leq r/2, t = 0\}$ and write $\overline{L} =
X + i Y$; then the vector fields $X$, $Y$, $[X,Y]$ span the tangent
space to $\p {\cal W}/S^1$ at each 
point of the complement of $J$, and are tangent to $J$ at every point
of $J$.  We conclude that the operator ${\cal L}_j$ on $\p {\cal
  W}/S^1$ is quite similar to the two-dimensional model that we
discussed above. 

There are two complications which we must note (and which are not
entirely trivial):  
{\bf (1)}  There are pseudodifferential factors, and the reduction of
the $\dbar$-Neumann problem to ${\cal L}$, and thereafter to ${\cal
  L}_j$, requires only a 
microlocal {\it a priori} estimate for ${\cal L}_j$ in a conic subset
of phase space; {\bf (2)}  The lower order terms $B_1 \overline{L}$,
$B_2 L, B_3$ are {\it not} negligible, indeed they determine  
the values of the exceptional Sobolev exponents, but the analysis can
be carried out for these terms as well.

It should be noted that a special feature of the worm is that the
rotational symmetry in $z_2$ makes possible (as we have noted) a
reduction to a 2-dimensional analysis, and this in turn produces a
certain convenient ellipticity. There is  no uniformity of estimates
with respect to $j$, but the analysis can be performed for each fixed
$j$. 

\bigskip

\section{The Automorphism Group of the Worm Domain}
\label{6}

It is of interest to know whether a biholomorphic mapping
of the smooth worm $\sW$ to any other smoothly bounded
pseudoconvex domain will extend to a diffeomorphism of
the closures.  Of course the
worm does {\it not} satisfy Condition $R$, so the obvious
tools for addressing this question are not available. As a
partial result, So-Chin Chen \cite{Che1} has
shown that the automorphism group of $\sW$ reduces to the
rotations in the $z_2$-variable; hence all biholomorphic {\it self-maps} of
$\sW$ do extend smoothly to the boundary.   
His result is this:
	      
\begin{thm}  \sl
Let $\sW$ be the smooth worm.  Then any automorphism
(i.e., biholomorphic self-map) of $\sW$ must be a rotation
in the $z_2$ variable.  In particular, the automorphism must extend
to a diffeomorphism of the closure.
\end{thm}
\proof
This is an interesting calculation.
First recall that, by Proposition \ref{pseudoconvexity}, 
the boundary of the smooth
worm $\sW$ is strongly pseudoconvex except on the annulus
$\mathcal A$ of the points
$(0,e^{i\log |z_2|^2})$ for $\big|\log |z_2|^2 \big|\le\mu$.

Let now $g = (g_1, g_2)$ be an automorphism of the worm 
$\sW$.  Then, by the fundamental result of Bell [Bel2], 
$g$ can be extended smoothly to all the strongly pseudoconvex points
of the  boundary.  In other words, $g$ extends to a
$C^\infty$-diffeomorphism of $\sW\setminus\mathcal A$ onto itself.
 
Consider now, for $e^{-\mu/2} < a < e^{\mu/2}$, the set
$$
T_a = \{(z_1,z_2) \in\p\sW: e^{-\mu/2} < |z_2| = a < e^{\mu/2}, z_1
\ne 0\} \, . 
$$
Notice that $T_a$ is a {\em deleted torus}, made of points of strong
pseudoconvexity in $\p\sW$.  Then, $g\bigr |_{T_a}$ is a
$C^\infty$-diffeomorphism having image contained in
$\p\sW\setminus\{(z_1,z_2)\in\p\sW:\, z_1\neq0,\,
\big|\log|z_2|^2\big|<\mu\}$.  

Then, if $\rho(z_1,z_2) = \big|z_1 - e^{i\log |z_2|^2} \big|^2 - (1 -
\eta(\log|z_2|^2))$  
is the obvious defining function for $\sW$, we see that the conditions
$$
|z_2|=a \qquad\text{and}\qquad
\biggl |g_1(z_1,z_2) - e^{i \log |g_2(z_1,z_2)|^2} \biggr | = 1
$$
define $T_a$. This implies that $\log|g_2(z_1,e^{i\theta}a)|^2 = \log
a^2$ for all $\theta$ real.
Thus $|g_2(z_1,z_2)|^2 = |z_2|^2 \cdot e^{2k\pi}$ for some
integer $k$.  By considering points $(z_1,z_2) \in T_a$ with $a$ close to
either $e^{-\mu/2}$ or $e^{\mu/2}$, we may conclude that $k = 0$ and
$|g_2(z_1,z_2)| = |z_2|$ for $(z_1,z_2) \in T_a$, 
$e^{-\mu/2}<a<e^{\mu/2}$.

Now, with $a$ fixed as before, let $\tilde z_1$ be a point with 
$|\tilde z_1 - e^{i \log |a|^2}| < 1$ and so that
$\tilde z_1$ lies in a small open neighborhood of $2e^{i\log |a|^2}$.

Consider the set given by
$$
 \{(\tilde z_1, z_2) \in \bC^2\} \cap \sW \, .
$$
It is not difficult to see that such a set is the union of (finitely
many)
concentric annuli.
Let $A_{\tilde z_1}$ denote one of these annuli, 
the inner boundary of $A_{\tilde z_1}$ being some circle $C_\alpha = 
\{(\tilde z_1, z_2) \in \p \sW:
|z_2| = \alpha\}$ and the outer boundary being another circle $C_\beta =
\{(\tilde z_1, z_2) \in \p \sW: 
|z_2| = \beta\}$ with $\alpha < a < \beta$.  Then $A_{\tilde z_1}$ is easily
identified with 
the planar annulus $A = \{z_2 \in \bC: \alpha < |z_2| < \beta\}$.  We
therefore obtain 
(using this identification)
\begin{equation}\label{eqno10}
g_2(\tilde z_1, C_\alpha) 
= C_\alpha \quad \hbox{and} \quad g_2(\tilde z_1, C_\beta)
= C_\beta \, .  
\end{equation}
Note that $g_2(\tilde z_1, \ \cdot \ )$ can be extended to an entire
function on all of $\bC$ just 
by using the Schwarz reflection principle.

Now line (\ref{eqno10}) tells us that
$$
g_2(z_1,z_2) = e^{i\theta(z_1)} \cdot z_2
$$
for some real function $\theta(z_1)$.  But $g_2$ is holomorphic in
$z_1$, we we may conclude that $\theta$ is a real constant $\theta_0$.
Thus 
\begin{equation}\label{eqno11}
g_2(z_1,z_2) = e^{i\theta_0} \cdot z_2 \quad \hbox{for} \ (z_1,z_2) \in \sW \, .
\end{equation}

As a consequence, since $g:\sW\raw\sW$ we have that
$$
\rho\bigl(g_1(z_1,z_2),g_2(z_1,z_2)\bigr)
= \big| g_1(z_1,z_2) - e^{i\log|g_2(z_1,z_2)|^2} \big|^2
< 1-\eta(\log|g_2(z_1,z_2)|^2) ;
$$
that is,
\begin{equation} \label{reduction}
 \big| g_1(z_1,z_2) - e^{i\log |z_2|^2} \big|^2
< 1-\eta(\log|z_2|^2) \ .
\end{equation}

Now examine the open, solid torus $\Pi_a$ given by
$$
\Pi_a = \{(z_1,z_2) \in \sW: e^{-\mu/2} < |z_2| = a < e^{\mu/2},\, 
|z_1 - e^{i\log |z_2|^2}| < 1\} \, . 
$$
Set, for $\theta$ real,
$$
\bigtriangleup_{a,\theta} 
= \{(z_1, a e^{i\theta}) \in \sW: |z_1 - e^{i\log  |a|^2}| < 1\} \, . 
$$
By (\ref{reduction})
it follows that the restriction of $g_1$ to $\bigtriangleup_{a,
  \theta_1}$ must map 
$\bigtriangleup_{a,\theta_1}$ biholomorphically onto 
$\bigtriangleup_{a,\theta_2}$  for some $\theta_2$.  
Thus the restriction of $g_1$ to $\bigtriangleup_{a,\theta_1}$
can be extended smoothly to $\overline{\bigtriangleup_{a,\theta_1}}$.  
We know that $g_1(0, a e^{i\theta_1}) = 0$, so it follows that $g_1(z_1,
z_2)$ can be expressed 
by way of the well-known automorphisms of the unit disc (see
\cite{GKr3}).  We see then that  
$$
g(z_1,ae^{i\theta})
   = e^{i\theta_0} \frac{b-\bigl(z_1-e^{i\log|a|^2}\bigr)}{1-\bar
  b\bigl( z_1 -e^{i\log|a|^2}\bigr)}  +e^{i\log|a|^2} 
$$
for some $\theta$ real and $b=b(ae^{i\log|a|^2})$ with $|b|<1$.  Using
the fact that
$g_1(0, a e^{i\theta_1}) = 0$, we calculate $\theta_0$ and we obtain that
\begin{equation}\label{eqno12}
g_1(z_1,z_2) = e^{i\log|z_2|^2} \left ( \frac{1 + \overline{b(z_2)} e^{i \log
      |z_2|^2}}{e^{i \log |z_2|^2} +b(z_2)} \right ) 
     \left ( \frac{b(z_2) + e^{i\log |z_2|^2} - z_1}{1 
- \overline{b(z_2)}(z_1 -
         e^{i\log |z_2|^2})} \right ) 
         - e^{i\log |z_2|^2} \, ;  
\end{equation}
here $b(z_2)$ is a real analytic function satisfying $|b(z_2)| < 1$
for $e^{-\mu/2} < |z_2| < e^{\mu/2}$. 
Equation (\ref{eqno12}) shows that there is a small $\varepsilon > 0$ such that 
$g_1(z_1,z_2)$ is real analytic on $\bigtriangleup(0,\varepsilon)
\times A_\d$, 
where $A_\d = \{z_2 \in \bC: e^{\d +\mu/2}
< |z_2| < e^{-\d +\mu/2}\}$, for
some small $\d > 0$. 
 Thus we see that $g_1(z_1,z_2)$ is holomorphic
on $\bigtriangleup(0,\varepsilon) \times A_\d$. 

As a consequence, one can write
$$
g_1(z_1,z_2) = \sum_{j=1}^\infty a_j(z_2) z_1^j \, ,
$$
with $a_j(z_2)$ holomorphic on $A_\d$ for all $j \geq 1$.  Direct
calculation yields that
$$
a_1(z_2) = \frac{\p g_1}{\p z_1} (0,z_2) = \frac{1 -
  |b(z_2)|^2}{|1 + \overline{b(z_2)} e^{i\log |z_2|^2}|^2} \, . 
$$
Thus $a_1(z_2)$ is a positive real constant, i.e., $a_1(z_2) = c > 0$.

Next we turn to the computation of $a_2(z_2)$.  Now
\begin{equation}\label{eqno13}
a_2(z_2) = \frac{1}{2} \cdot \frac{\p^2 g_1}{\p z_1^2}(0,z_2) =
c \cdot \frac{\overline{b(z_2)}}{1 + \overline{b(z_2)}  
     e^{i\log |z_2|^2}} \, .  
\end{equation}
We assert that $a_2 \equiv 0$.	Now set
$$
h(z_2) = \frac{a_2(z_2)}{c} = \frac{\overline{b(z_2)}}{1 +
  \overline{b(z_2)}e^{i\log|z_2|^2}} \, , 
$$
which is holomorphic on $A_\d$.  We see that
\begin{eqnarray*}
c & = & 
\frac{1 - |b(z_2)|^2}{|1 + \overline{b(z_2)}e^{i\log |z_2|^2}|^2} \\
  & = & |1 + h(z_2) e^{i \log |z_2|^2}| - |h(z_2)|^2 \\
  & = & 1 + 2 \hbox{Re} (h(z_2) e^{i \log |z_2|^2}) \, .
\end{eqnarray*}

In conclusion, we may write
$$
h(z_2) e^{i \log |z_2|^2} = c_0 + i I(z_2) \, ,
$$
where $c_0 = \frac{1}{2} (c-1)$ and $I(z_2)$ is a smooth, real-valued
function on $A_\d$. 
Thus we have
\begin{equation}\label{eqno14}
h(z_2) = c_0 e^{-i \log |z_2|^2} + i I(z_2) e^{-i \log |z_2|^2} \, ,
\end{equation}
which is holomorphic on $A_\d$.  Locally we may multiply equation
(\ref{eqno14}) by $e^{2i\log z_2}$  
to obtain a new holomorphic function, and we find that
$$
h(z_2) e^{2i \log z_2} = c_0 
e^{-2{\rm arg} \, z_2} + i I(z_2) e^{-2 {\rm arg} \, z_2} \, .
$$

Of course the real part of $g(z_2) e^{2i \log z_2}$ is a harmonic function.
We write 
$z_2 = u + iv$ as usual.  By direct computation we find that
\begin{eqnarray*}
\d_{z_2} (c_0 e^{-2 {\rm arg} \, z_2}) & = & c_0
\d_{z_2}(e^{-2\tan^{-1} v/u}) \\ 
    & = & \frac{4c_0}{u^2 + v^2} e^{-2\tan^{-1} v/u} \equiv 0 \, .
\end{eqnarray*}
This entails $c_0 = 0$ so that $c = 1$.  Thus (\ref{eqno14}) reduces to  
$$
-i h(z_2) = I(z_2) e^{-i \log |z_2|^2} \, ,
$$ 
which is holomorphic on $A_\d$.  Repeating the very same argument,
we find that 
$$
- i h(z_2) e^{2i \log z_2} = I(z_2) e^{-2 {\rm arg} \, z_2} = c_1 \, ,
$$
where $c_1$ is a constant.  Thus
$$
I(z_2) = c_1 e^{2{\rm arg}\, z_2} 
$$
is a well-defined function on $A_\d$.  This forces $c_1 = 0$.  As
a result, 
$h(z_2) \equiv 0$.  We see in sum that $a_2(z_2) = 0$ as claimed.

Now we may conclude from (\ref{eqno13}) that $b(z_2) \equiv 0$ on $A_\d$.
Therefore 
equation (\ref{eqno12}) simplifies to
\begin{equation}\label{eqno15}
g_1(z_1,z_2) \equiv z_1 \quad \hbox{on} \  \sW \, . 
\end{equation}
The result now follows from (\ref{eqno11}) and (\ref{eqno15}).
\endpf

\bigskip 

\section{Analysis on $D_\beta$ and $D'_\beta$}\label{7}

We now summarize our main results about the non-smooth
worm domains $D_\beta$ and $D'_\beta$.  Details
appear in \cite{KrPe}.
		
It is of interest, in its own right and as a model
for the smooth case, to study the behavior of the Bergman
kernel and projection on the non-smooth worm domain $D_\b$ and its
biholomorphic copy $D_\b'$. 
An important transformation rule for the Bergman kernel
(see \cite{Kr1}) says that if $\Phi: \Omega_1 \raw \Omega_2$ is
biholomorphic then
$$
K_{\Omega_1}(z, \zeta) = \hbox{det} \, \hbox{Jac}_{\Phi_1}(z) \cdot
K_{\Omega_2}(\Phi(z), \Phi(\zeta)) 
       \cdot \overline{\hbox{det} \, \hbox{Jac}_{\Phi_1}(\zeta)} \, .
$$
It is obvious from this transformation rule 
that it suffices to obtain the expression for the kernel in
just one of the two domains (see (\ref{eq:biholo})). 
However, 
the $L^p$-mapping properties
of the Bergman projections of the two domains turn out to be
substantially different (just because $L^p$ spaces of
holomorphic functions do not transform canonically under
biholomorphic maps when $p \ne 2$, due to the presence of the Jacobian
factor).  We shall explore this
result in the next Theorems \ref{boundedness-on-Dbeta-prime} 
and \ref{boundedness-on-Dbeta}.

We shall discuss here the results (contained in detail in
\cite{KrPe}) 
concerning the explicit expression of the Bergman
kernels for $D_\b$ and $D_\b'$.  Once these are available we study the 
$L^p$-mapping properties of the corresponding Bergman
projections.  More precisely we prove the following theorems. 
There are two principal results of \cite{KrPe}.  Of course
the proofs, which are quite technical, must be omitted.  But the
applications that we provide give a sense of the meaning and
significance of these two theorems.

\begin{thm}\label{thm1}	 \sl
Let $c_0$ be a positive fixed constant.  Let
$\chi_1$ be a smooth cut-off function on the real line,
supported on $\{x:\, |x|\le 2c_0\}$,
identically 1 for $|x|<c_0$. Set $\chi_2 =1-\chi_1$.

Let $\b>\pi$ and let $\nu_\b = \pi/[2\beta - \pi]$. 
Let $h$ be fixed, with
\begin{equation}\label{(*)}
\nu_\b<h<\min(1,2\nu_\b) \, .  
\end{equation}
Then there exist
functions $F_1,F_2,\dots, F_8$ 
and $\tilde F_1,\tilde F_2,\dots, \tilde F_8$,
holomorphic in $z$ and anti-holomorphic in $w$, for $z=(z_1,z_2)$,
$w=(w_1,w_2)$ varying in a
neighborhood of $D'_\b$, and 
having size ${\cal O}(|\Re z_1-\Re w_1|)$, together with all their
derivatives, for $z,w \in \overline{D'_\b}$,
as $|\Re z_1-\Re w_1|\raw+\infty$.  
Moreover, there exist functions $E,\tilde E \in {\mathcal C}^\infty 
(\overline{D'_\b}\times \overline{D'_\b})$ such that
$$
D_{z_1}^\alpha D_{w_1}^\gamma E(z,w),
\, D_{z_1}^\alpha D_{w_1}^\gamma \tilde E(z,w)
=\mathcal O(|\Re z_1-\Re w_1|^{|\alpha|+|\gamma|})\ ,
$$
as $|\Re z_1-\Re w_1|\raw+\infty$. (Here, for $\lambda\in\bC$,
$D_\lambda$ denotes the 
partial derivative in $\l$ or $\bar\l$.)

Then the following holds.  Set
\begin{align}
K_b(z,w)
& =   \frac{F_1(z,w)
}{(i(z_1- \overline{w}_1)+2\b)^2 (e^{(\beta - \pi/2)} - z_2 \overline{w}_2)^2}  
\notag \notag\\    
&\qquad\quad
+ \frac{F_2(z,w)}{(i(z_1- \overline{w}_1)+2\b)^2
(z_2 \overline{w}_2 - e^{-[i(z_1-  \overline{w}_1) +\pi]/2})^2 }
\notag \notag\\  
& \quad
+
\frac{F_3(z,w)}{(e^{[\pi-i(z_1- \overline{w}_1)]/2}-z_2 \overline{w}_2)^2
(e^{(\beta - \pi/2)}-z_2 \overline{w}_2)^2}
\notag \notag\\
& \qquad +  
\frac{F_4(z,w)}{(i(z_1- \overline{w}_1)-2\b)^2 
(e^{[\pi - i(z_1- \overline{w}_1)]/2} - z_2 \overline{w}_2)^2}
\notag \notag\\
&\qquad\quad
+ \frac{F_5(z,w)}{(i(z_1- \overline{w}_1)-2\b)^2  
(z_2 \overline{w}_2 - e^{-(\beta - \pi/2)})^2}  
\notag \notag\\
&\quad 
+ \frac{F_6(z,w)
}{(e^{-[i(z_1- \overline{w}_1)+\pi]/2}-z_2 \overline{w}_2)^2
(e^{-(\beta - \pi/2)}-z_2 \overline{w}_2)^2} 
 \notag \notag\\
& \quad
+ \frac{F_7(z,w)}{
(i(z_1- \overline{w}_1)+2\b)^2(e^{(\beta - \pi/2)} - z_2 \overline{w}_2)
(e^{-[i(z_1- \overline{w}_1) + \pi]/2}-z_2 \overline{w}_2)} 
\notag \notag\\
& \quad
+ \frac{F_8(z,w)}
{(i(z_1- \overline{w}_1)-2\b)^2  
(e^{-(\beta - \pi/2)}-z_2 \overline{w}_2)(e^{[\pi-i(z_1-
  \overline{w}_1)]/2}-z_2 
\overline{w}_2)}  
+E(z,w) \notag\\ 
& \equiv K_1(z,w)+\cdots+K_8(z,w) + E(z,w)\ .
\label{K1throughK8}
\end{align}
Define $K_{\tilde b}$ by replacing $F_1,\dots,F_8$ and $E$ by
$\tilde F_1,\dots,\tilde F_8$ and $\tilde E$ 
and thus $K_1,\dots,K_8$  by
$\tilde K_1,\dots,\tilde K_8$ respectively in formula
(\ref{K1throughK8}).  

Then there exist functions $\phi_1,\phi_2$ entire in $z$ and
$\overline{w}$ (that is, anti-holomorphic in 
$w$), which are of size ${\cal O}(|\Re z_1 -\Re w_1|)$, together with
all their 
derivatives, uniformly in all closed strips 
$\{|\Im z_1|+|\Im w_1|\le C\}$, such that  
the Bergman kernel $K_{D'_\b}$ on $D'_\b$ admits the asymptotic expansion
\begin{multline}
K_{D'_\b}(z,w) =
\chi_1(\Re z_1 -\Re w_1) K_b(z,w) + 
\chi_2(\Re z_1 -\Re w_1)
\biggl\{ 
e^{-h\sgn (\Re z_1 -\Re w_1) \cdot (z_1- \overline{w}_1)}   
K_{\tilde b}(z,w) 
\\
+e^{-\nb \sgn (\Re z_1 -\Re w_1)\cdot  (z_1- \overline{w}_1)}
\biggl( 
 \frac{\phi_1 (z_1,w_1)}{(e^{[\pi- i(z_1- \overline{w}_1)]/2}- z_2 \overline{w}_2)^2} 
   + \frac{\phi_2 (z,w)}{
(e^{-[i(z_1- \overline{w}_1)+\pi]/2}-z_2 \overline{w}_2)^2} \biggr) \biggr\} \, .
\label{K-Dbetaprime}
\end{multline}
Here $h$ is specified as in (\ref{(*)})  above.
\end{thm}

\begin{thm}\label{thm2}	 \sl
With the notation as in Theorem \ref{thm1}, 
there exist
functions $g_1,g_2$, 
$G_1,G_2,\dots, G_8$ and $\tilde G_1,\tilde G_2,\dots, \tilde G_8$,
holomorphic in $\z$ and anti-holomorphic in $\om$, for $\z=(\z_1,\z_2)$,
$\om=(\om_1,\om_2)$ varying in  $\overline{D'_\b}\setminus\{(0,z_2)\}$, 
such that
$$
\p_{\z_1}^\alpha \p_{\overline{\om}_1}^\gamma G(\z,\om)=
{\cal O} \bigl( |\z_1|^{-|\alpha|} |\om_1|^{-|\gamma|} \bigr)
\qquad\qquad\text{as}\quad |\z_1|,|\om_1|\raw 0\ ,
$$
where $G$ denotes any of the functions $g_j$, $G_j$, $\tilde G_j$.
Moreover, there exist functions $E,\tilde E \in {\mathcal C}^\infty
\bigl( \overline{D'_\b}\setminus\{(0,z_2)\} \times
\overline{D'_\b}\setminus\{(0,z_2)\} \bigr)$ such that
$$
D_{\z_1}^\alpha D_{\om_1}^\gamma E(\z,\om),
D_{\z_1}^\alpha D_{\om_1}^\gamma \tilde E(\z,\om) =
{\cal O} \bigl( |\z_1|^{-|\alpha|} |\om_1|^{-|\gamma|} \bigr)
\qquad\qquad\text{as}\quad |\z_1|,|\om_1|\raw 0\ .
$$
(Here $D_\lambda$, 
for $\lambda\in\bC$,
$D_\lambda$ denotes the 
partial derivative in $\l$ or $\bar\l$.)

Then the following holds.  Set
\begin{align} 
H_b(\z,w)
& =   \frac{G_1(\z,w)
}{(i\log(\z_1/\overline{\om}_1)+2\b)^2 (e^{(\beta - \pi/2)} - \z_2
  \overline{\om}_2)^2}   
\notag\\    
&\qquad\quad
+ \frac{G_2(\z,w)}{(i\log(\z_1/\overline{\om}_1)+2\b)^2
\bigl((\z_1/\overline{\om}_1)^{-i/2}e^{-\pi/2} 
-\z_2 \overline{\om}_2  
\bigr)^2 }
\notag\\  
& \quad
+
\frac{G_3(\z,w)}{\bigl((\z_1/\overline{\om}_1)^{-i/2}e^{\pi/2} 
-\z_2 \overline{\om}_2\bigr)^2
(e^{(\beta - \pi/2)}-\z_2 \overline{\om}_2)^2}
\notag\\
& \qquad +  
\frac{G_4(\z,w)}{(i\log(\z_1/\overline{\om}_1)-2\b)^2 
\bigl((\z_1/\overline{\om}_1)^{-i/2}e^{\pi/2} 
- \z_2 \overline{\om}_2\bigr)^2}
\notag\\
&\qquad\quad
+ \frac{G_5(\z,w)}{(i\log(\z_1/\overline{\om}_1)-2\b)^2  
(e^{-(\beta - \pi/2)}-\z_2 \overline{\om}_2)^2}  
\notag\\
&\quad 
+\frac{G_6(\z,w)
}{\bigl((\z_1/\overline{\om}_1)^{-i/2}e^{-\pi/2} 
-\z_2 \overline{\om}_2\bigr)^2
(e^{-(\beta - \pi/2)}-\z_2 \overline{\om}_2)^2} 
 \notag\\
&\quad
+ \frac{G_7(\z,w)}{ 
(i\log(\z_1/\overline{\om}_1)+2\b)^2(e^{(\beta - \pi/2)} - \z_2
\overline{\om}_2) 
\bigl((\z_1/\overline{\om}_1)^{-i/2}e^{-\pi/2}
-\z_2 \overline{\om}_2\bigr)} 
\notag\\
&\quad
+ \frac{G_8(\z,w)}
{(i\log(\z_1/\overline{\om}_1)-2\b)^2  
(e^{-(\beta - \pi/2)}-\z_2
\overline{\om}_2)\bigl((\z_1/\overline{\om}_1)^{-i/2}e^{\pi/2}  
-\z_2 \overline{\om}_2\bigr)}  
+E(\z,w)\notag\\
& \equiv H_1(\z,\om)+\cdots+H_8(\z,\om) + E(\z,\om)\ . 
\end{align}
Define $H_{\tilde b}$ by replacing $G_1,\dots,G_8$ and $E$ by
$\tilde G_1,\dots,\tilde G_8$ and $\tilde E$, and $H_1,\dots,H_8$  by
$\tilde H_1,\dots,\tilde H_8$,  respectively.

Then, setting $t=|\z_1|-|\om_1|$, we have this
asymptotic expansion for the Bergman kernel on $D_\b$:
\begin{align*}
& K_{D_\b} \bigl((\z_1,\z_2),(\om_1,\om_2)\bigr)\\
&= \chi_1(t)
\frac{H_b(\z,\om)
}{\z_1 \overline{\om}_1} +
\chi_2(t) 
\biggl\{  \biggl(\frac{|\z_1|}{|\om_1|}\biggr)^{-h\sgn t} 
e^{-h\sgn t \cdot (\arg\z_1+\arg\om_1)}   
\frac{H_{\tilde b}(\z,\om)}{ \z_1
  \overline{\om}_1}  
  \notag\\
& \quad
+ \biggl(\frac{|\z_1|}{|\om_1|}\biggr)^{-\nu_\b\sgn t} 
e^{-\nu_b\sgn t \cdot (\arg\z_1+\arg\om_1)}\biggl(    
\frac{g_1 (\z_1,\om_1)}{ \z_1\overline{\om}_1}\cdot
\frac{1}{\bigl(
( \z_1/\overline{\om_1})^{-i/2}
e^{\pi/2}- z_2 \overline{\om}_2 \bigr)^2} \notag\\
& \quad\qquad\qquad   
+\frac{ g_2 (\z,\om)}
{\z_1\overline{\om}_1} \cdot \frac{1}{\bigl( 
( \z_1/\overline{\om}_1)^{-i/2}
e^{-\pi/2}- z_2 \overline{\om}_2 \bigr)^2} \biggr)  \biggr\}
\, ,
\end{align*}
where $h$ is defined in (\ref{(*)}).
\end{thm}

The Bergman projection is trivially
bounded on $L^2(\Omega)$.  It is of some interest to ask about
the mapping properties of $P$ on $L^p(\Omega)$ for $1 \leq p \leq \infty$.
In general
$P$ will be bounded on either $L^1$ or $L^\infty$. 
This assertion follows because 
$P$ is in the nature of a Hilbert integral,
see \cite{PhS1}, \cite{PhS2}.  Details are provided in that source.

Matters for $L^p$, with $1<p<\infty$,  are more subtle.
In fact the results are different
for the two domains $D'_\beta$ and $D_\beta$.  We will use our
asymptotic expansions to prove the following theorems:

\begin{thm}\label{boundedness-on-Dbeta-prime}   \sl
The Bergman projection $P_{D'_\beta}$ on the domain $D'_\beta$ is bounded on
$L^p$ for $1 < p < \infty$.
\end{thm}

\begin{thm}\label{boundedness-on-Dbeta}    \sl
Let $\beta > \pi$ and $\nu_\beta = \pi/[2\beta - \pi]$.  
The Bergman projection $P_{D_\beta}$ on the domain $D_\beta$ 
is only bounded on $L^p$ for $2/[1 + \nu_\b] \leq p \leq 2/[1 - \nu_\b]$.  It
is unbounded on $L^p$ for $p$ outside this range.
\end{thm}
						 
This situation is at first puzzling because the two domains
$D'_\beta$ and $D_\beta$ are biholomorphic. But, whereas it is
well known that $L^2$ transforms canonically under
biholomorphic maps (see \cite{Kr1}), such is not the case for
$L^p$ when $p \ne 2$. 
						 
We shall now describe the proof of the result on $D_\beta$, which
is the most interesting case.  
In fact, we concentrate on the negative part of the result, since it
bears some consequences on the unboundedness of the Bergman projection
of $\mathcal W$.  The proof positive part is more direct, and it
relies on an involved, systematic application of Schur's Lemma.

The
proof of the cognate result for $D'_\beta$ is 
somewhat more elementary, although 
too elaborate to
present here.
\medskip

We concentrate the on the negative result.
Let $1 < p < \infty$ and assume that $P:L^p(D_\b) \rightarrow
L^p(D_\b)$ is bounded. It follows that
for any $\z\in D_\b$ fixed, 
$K_{D_\b} (\cdot,\z) \in L^{p'}(D_\b)$, where
$p' = p/[p-1]$ is the exponent conjugate to $p$.  

For, if $P=P_{D_\b}$ is bounded, then
for all $f\in L^p(D_{\b})$ and all $\z\in D_\b$,
$$ 
|\langle f, K_{D_\b}(\cdot, \z) \rangle | =  |Pf(\z)| 
\leq  c_\z \|Pf\|_{L^p} \le C \|f\|_{L^p}\, .
$$

\begin{lemma}\label{unboundedness}	\sl
For any $\z\in D_\b$ it holds that
$K_{D_\b}(\cdot,\z)\in L^p(D_\b)$ only if
$2/(1 + \nu_\b) < p < 2/(1 - \nu_\b)$.
\end{lemma}
\proof
Fix $\z\in D_\b$ and define
$$
\Omega_\z = \bigl\{\om \in D_\b: |\om_1| < |\z_1| \ , \
1/4 \leq | e^{\pi/2} (\z_1/\overline{\om}_1 )^{\pm i/2} - 
\z_2 \overline{\om}_2 |
          \leq 1/2\bigr\} \, .
$$
Recall the expansion for the kernel $K_{D_\b}(\z,\om)$ given in
Theorem \ref{thm2}.  Then, for $\om\in\Omega_\z$,
we have that 
$$
|H_b(\z,\om)|,\, |H_{\tilde{b}}(\z,\om)|\le C_\z
$$
for some constant independent of $\om$, so that
\begin{equation}\label{estimate-from-below-K-D-beta}
|K_{D_\b}(\z,\om)|
\ge c_\z |\om_1|^{\nu_\b -1}
\end{equation}
for $\om\in\Omega_\z$.  

Therefore
\begin{align*}
\int_{D_{\b}} |K_{D_\b} (\cdot,\z)|^{p'} \, dV(\om) 
& \ge 
    \int_{\Omega_\z} |K_{D_\b} (\cdot,\z)|^{p'} \, dV(\om)\\
& \ge  c_\z 
\int_{\Omega_\z} \bigl(|\om_1|^{\nu_\b-1} 
\bigr)^{p'} \, dV(\om) 
 = c \int_0^{|\z_1|} r^{p'(\nu_\b-1) +1} \, dr \ .
\end{align*}

Obviously for convergence we need $p' (\nu_\b -1)+1>-1 $, that is
$p' < 2/[1 -\nu_\b]$.  Hence if
$p \geq 2/[1 - \nu_\b]$ then the integral diverges.  The other
result, for $p \leq 2/[1 + \nu_\b]$, follows by duality.
This proves the lemma.
\endpf

We now show 
that we can use
Barrett's exhaustion procedure (see \cite{Ba2}) to obtain a
negative result with the same indices on the smooth worm
$\mathcal W$. 

\begin{thm}\label{straube-thm}	 \sl
Let $\mathcal P$ denote the Bergman projection on the smooth, bounded
worm $\mathcal W=\mathcal W_\b$, with $\b>\pi/2$.  Then, if $\mathcal P:
L^p(\mathcal W_\b)\raw L^p(\mathcal W_\b)$ is bounded, necessarily
$2/[1+\nu_\b]<p<2/[1-\nu_\b]$. 
\end{thm}

\proof
Suppose 
$\mathcal P:
L^p(\mathcal W)\raw L^p(\mathcal W)$ is bounded for a given
$p\neq2$.  Without loss of generality, we may assume that $p>2$.

Let $\tau_R$ be defined as $\tau_R(z_1,z_2)=(Rz_1,z_2)$, for
$R\ge1$. Recall that 
$$
\mathcal W' =\bigl\{
(z_1,z_2):\ \big|z_1 -e^{i\log|z_2|^2}\big|<1,\ 
\big|\log|z_2|\big|<\beta-\pi/2 \, \bigr\}
$$
 is the truncated version of $\mathcal W$.  Then,
$$
\tau_R(\mathcal W)\supseteq 
\tau_R(\mathcal W') = 
\bigl\{
(w_1,w_2):\ |w_1|^2/R  - 2\Re \bigl(w_1 e^{-i\log|w_2|^2}\bigr) |<0,\ 
\big|\log|w_2|\big|<\beta-\pi/2 \, \bigr\}
$$
and 
$$
\tau_R(\mathcal W') \nearrow D_\beta\qquad\text{as}\quad R\raw+\infty\, .
$$

Let $T_R\, :L^p \bigl(\tau_R(\mathcal
W)\bigr) \raw L^p (\mathcal W)$ be defined as $T_R(f)=
f\circ\tau_R$. Notice that 
$$
\|T_R(f)\|_{L^p (\mathcal W)}^p
 =  R^{-2} 
\|f\|_{L^p (\tau_R(\mathcal W))}^p
$$
for all $f\in L^p (\tau_R(\mathcal W))$.
Moreover, set $P_R= T_R^{-1 }\mathcal P T_R$ and notice that $P_R$ is
the Bergman projection of $\tau_R(\mathcal W')$.
Using the boundedness
of $\mathcal P$ on $L^p(\mathcal W)$ it follows that
\begin{equation}\label{uniform-bound}
\| P_R \varphi\|_{L^p(\tau_R(\mathcal W'))}
\le C \|\varphi\|_{L^p(\tau_R(\mathcal W'))}
\end{equation}
for all continuous functions $\varphi$ 
with compact support in $\tau_R(\mathcal W')$, 
where $C$ is a constant independent
of $R$. \medskip

We now claim that $(P_R\varphi)\,\widetilde{\,}\raw
(P_{D_\beta}\varphi)\,\widetilde{\,}$ weakly in $L^p(\bC^2)$, as
$R\raw+\infty$. Here,  by $f\,\widetilde{\,}$ we denote the extension of
$f$ to all of $\bC^2$ (defining
$f\,\widetilde{\,}=0$  outside the natural domain of definition of $f$). 
Notice that, since 
$\tau_R(\mathcal W') \nearrow D_\beta$ as $R\raw+\infty$, if
$\varphi$ has compact support in 
$D_\beta$ then there exists $R_0$ such that $\varphi$ has compact support
in $\tau_R(\mathcal W')$, for $R\ge R_0$.\medskip

Assume the claim for now, and we finish the proof. 
For all $\varphi$ continuous
with compact support in $D_\beta$,
by
(\ref{uniform-bound})
 and the claim it follows at once that
$$
\|P_{D_\beta} \varphi\|_{L^p(D_\beta)} 
\le C \| \varphi\|_{L^p(D_\beta)} \ .
$$
The result now follows from Theorem \ref{boundedness-on-Dbeta}.
\medskip

Finally, we prove the claim.  Since
$\|(P_R\varphi)\,\widetilde{\,}\|_{L^p(\bC^2)} 
\le C \|\varphi\,\widetilde{\,}\|_{L^p(\bC^2)}$, there exists a
subsequence of $\{ (P_R\varphi)\,\widetilde{\,}\}$ that converges
weakly to $g\in L^p(\bC^2)$.  Notice that $g$ vanishes out
$\overline{D_\beta}$, since $(P_R\varphi)\,\widetilde{\,}$ does for
all $R\ge1$.  Moreover, $g$ is holomorphic on $D_\beta$, since 
$\tau_R(\mathcal W') \nearrow D_\beta$ 
and, $(P_{R_2}\varphi)\,\widetilde{\,}=
(P_{R_1}\varphi)\,\widetilde{\,}$  
on $\tau_{R_1}(\mathcal W')$
when $R_2>R_1$.
Moreover notice that, since 
$\|(P_R\varphi)\,\widetilde{\,}\|_{L^2(\bC^2)} 
\le \|\varphi\,\widetilde{\,}\|_{L^2(\bC^2)}$, there exists a
subsequence $\{ (P_{R_j} \varphi)\,\widetilde{\,}\}$ that converges weakly
also in $L^2(\bC^2)$, to the same function $g$.  We wish to show
that $g=P_{D_\beta}$ on $D_\beta$. 

Now let $h\in A^{p'}\cap A^2 (D_\beta)$ (where
$A^p$ denotes the Bergman space, and $p'$ is the exponent conjugate to
$p$).  We have that
\begin{align*}
\la g-P_{D_\beta}\varphi,\, h\ra_{D_\beta} 
& = \lim_{j\raw+\infty}  \la (P_{R_j}\varphi)\,\widetilde{\,}
-P_{D_\beta}\varphi,\, h\ra_{D_\beta}  \\
& = \lim_{j\raw+\infty}  \int_{\tau_{R_j}(\mathcal W')} 
 P_{R_j}\varphi \overline{h} 
-\int_{D_\beta} P_{D_\beta} \varphi \overline{h}   \\
& = \lim_{j\raw+\infty}  \int_{D_\beta\setminus \tau_{R_j}(\mathcal W')} 
\varphi \overline{h}  = 0\ .
\end{align*}
If we show that $A^p\cap A^2 (D_\beta)$ is dense in $A^2(D_\beta)$, it
would follow that $g-P_{D_\beta}\varphi\perp A^2(D_\beta)$, that is
$g=P_{D_\beta}$, and we would be done.

Notice that for $\delta>0$ we have 
\begin{align*} 
\iint_{D_\beta} \big|e^{-\delta z_1^2}\big|^q \, dz_1 dz_2 
& = \iint_{D_\beta'} \big|e^{-\delta \log \z_1^2}\big|^q \,
\frac{1}{|\z_1|^2}
d\z_1 d\z_2  \\
& \le c \int_{|\log |\z_2^2|<\mu} \int_{|\Im\z_1-\log|\z_2|^2|<\pi/2} 
e^{-\delta q\log^2|\z_1|}  \frac{1}{|\z_1|^2} d\z_1 d\z_2  \\
& \le c \int_0^{+\infty} 
e^{-\delta q\log^2r}  r^{-2} dr <\infty\ ,
\end{align*}
that is $e^{-\delta z_1^2}\in L^q(D_\beta)$ for all $q$'s.
Now, for any $h\in A^2(D_\beta)$, consider $h_\delta=
he^{-\delta z_1^2}$. 
Since
$e^{-\delta z_1^2}$ is bounded, 
$h_\delta\in A^2(D_\beta)$.  Moreover, 
since $p$ is taken to be larger than $2$, by H\"older inequality we have 
$\| h_\delta \|_{L^{p'}}\le c\| h\|_{L^2}^{1/p'} <\infty$. 
Thus $h_\delta\in A^2\cap A^{p'}(D_\beta)$ and clearly $h_\delta$
converges to $h$ in the $L^2$-norm.  
This completes the proof of the theorem.
\endpf
\vspace*{.12in}

\bigskip

\section{Irregularity Properties of the Bergman Kernel}\label{8}

We now examine the boundary asymptotics for the Bergman kernel on the
domains $D_\beta$ and $D_\b'$
and determine various irregularity properties of the corresponding
Bergman kernels.

Begin with the asymptotic formula in Theorem \ref{thm1}.
We point out 
particularly 
that there are two kinds of behavior: one kind at the "finite portion
of the boundary" and the other one as $|\Re z_1 - \Re _1|
\raw+\infty$.  These two different behaviors are expressed by the
first and second terms in (\ref{K-Dbetaprime}), respectively.
For the former type, we notice from (\ref{K1throughK8}) that 
the lead terms have expressions in the
denominator of products of two terms like
$$
(i(z_1\pm \overline{w}_1)+2\b)^2\, ,\quad 
(z_2 \overline{w}_2 -e^{\pm(\beta - \pi/2)})^2\, ,\quad\text{and}\quad   
(z_2 \overline{w}_2 - e^{-[i(z_1-  \overline{w}_1) \pm\pi]/2})^2 \, .
$$
These singularities are similar to the ones of a Bergman kernel of a
domain in $\bC^2$ which is essentially a product domain.  It is
important to observe that the kernel does not become singular only
when $z,w$ tend to the same point on the boundary.  For instance, 
it becomes singular as $(i(z_1\pm \overline{w}_1)+2\b)\raw0$, while
there is no restriction  on the behavior of $z_2$ and $w_2$.
We will be more detailed below in the case of the domain $D_\b$.
For the case of this domain, we finally notice that
the main term at infinity, that is when 
$|\Re z_1 - \Re _1|\raw+\infty$, behaves like
$e^{-\nu_\beta|z_1-\bar w_1|}\cdot(z_2\bar w_2)^{-1}$.\medskip

Next we consider the case of $D_\b$.  The mapping $(z_1,z_2)\in D_\b'
\mapsto (\z_1,\z_2)=(e^{z_1},z_2)\in D_\b$ sends the point at infinity
(in $z_1$) 
into the origin (in $\z_1$).  Keeping into account the Jacobian
factor, when $|\z_1|-|\om_1|\raw0^+$, 
the kernel on $D_\b$ is asymptotic to
$$
\frac{|\om_1|^{\nu_\b-1}}{|\z_1|^{\nu_\b+1}}\cdot (\z_2\bar\om_2)^{-1}\
.
$$

Recall 
the inequalities that define $D_\beta$:
$$
D_\b = \Big\{\z \in \bC^2: 
\Re \bigl(\z_1 e^{-i\log |\z_2|^2}
\bigr)>0 ,\, \big|\log |\z_2|^2\big| < 
\b - \frac{\pi}{2} \Big\} \, .
$$
If we take $\z,\,\om\in D_\b$ and
let $\om_1$ tend to 0, 
then clearly $\om\raw \p D_\b$ and
$\z_1,\,\z_2,\,\om_2$ are unrestricted. Therefore,
$K_{D_\b}(\cdot,\om) \not\in C^\infty(D_\b)$
for $\om\in\{ (0,\om_2)\}$, with $|\log|\om_2|^2|<\b-\pi/2$.

Notice that this is in contrast, for instance, to the
situation on the ball or, more generally, on a strongly pseudoconvex
domain.  On either of those types of domains $\Omega$, the kernel is known
to be smooth on $\overline{\Omega} \times \overline{\Omega} \setminus 
(\Delta \cap [\p \Omega \times \p \Omega])$.  See
\cite{Ker} and \cite{Kr1}.

By the same token (by almost the same calculation), it is easy to 
conclude that the Bergman
projection on 
$D_\beta$ {\it cannot} map functions in $C^\infty(\overline{D}_\beta)$ to
functions in $C^\infty(\overline{D}_\beta)$.   This, of course, is the
failure of Condition $R$ on these domains.\medskip

In Section \ref{ConditionR} we have seen that
$P_{\sW}:C^\infty(\overline{\sW}) \not\raw
C^\infty(\overline{\sW})$, that is, that $\sW$ does not
satisfy Condition $R$. A philosophically related fact, due to
Chen \cite{Che2} and Ligocka \cite{Lig} independently, is that
the Bergman kernel of $\sW$ cannot lie in
$C^\infty(\overline{\sW}\times \overline{\sW}\setminus\Delta)$
(where $\Delta$ is the boundary diagonal). In fact, in \cite{Che2}
it is shown that this phenomenon is a consequence of the
presence of a complex variety in the boundary of $\sW$. We
shall also explore this singularity phenomenon in what follows. 
\medskip

The proof of the general result of So-Chin Chen
follows a classical paradigm for establishing propagation of singularities
for the $\dbar$-Neumann problem and similar phenomena.  

\begin{thm}    	  \sl
Let $\Omega \subseteq \bC^n$ be a smoothly bounded, pseudoconvex
domain with $n \geq 2$.  Assume that there is a complex variety $V$,
of complex dimension at 
least 1, in $\p \Omega$.   Then
$$
K_\Omega(z,w) \not \in C^\infty(\overline{\Omega} \times
\overline{\Omega} \setminus \bigtriangleup(\p \Omega)) \, , 
$$
where $\bigtriangleup(\p \Omega) = \{(z,z): z \in \p \Omega\}$.  
\end{thm}
\proof
Let $p \in V$ be a regular point.  Let $n_p$ be the unit outward
 normal vector at $p$.  Then there are small numbers $\d,
 \ve_0 > 0$ such that 
 $w - \ve n_p \in \Omega$ for all $w \in \p \Omega \cap
 B(p, \d)$ and 
 all $0 < \ve < \ve_0$.  Let {\bf d} be an analytic disc in
 $\p \Omega \cap B(p,\d) \cap V$. 
We may assume that this disc is centered at $p$.  In other words, {\bf
  d} is the image of  
the unit disc in the plane mapped into $\bC^n$ with the origin going
to $p$. 

Seeking a contradiction, we assume that $K_\Omega(z,w) \in
C^\infty(\overline{\Omega} \times \overline{\Omega} \setminus  
\bigtriangleup(\p \Omega))$.  Then we certainly have
$$
\sup_{w \in \p {\bf d}} |K_\Omega(p,w)| \leq M < + \infty   \eqno (22)
$$
for some positive, finite number $M$.  On the other hand, we know
(see \cite{BlP}) that
$$
\lim_{\ve \raw 0} K_\Omega (p - \ve n_p, p - \ve
n_p) = + \infty \, .  \eqno (23) 
$$
By the maximum modulus principle we then obtain
$$
\sup_{w \in \p {\bf d}_\ve} |K_\Omega(p - \ve n_p, w)|  \geq
K_\Omega (p - \ve n_p, p - \ve n_p) \, , 
$$
where ${\bf d}_\ve = {\bf d} - \ve n_p \subseteq \Omega$.
We conclude that 
$$
\sup_{w \in \p {\bf d}} |K_\Omega (p,w)| = \lim_{\ve \raw 0}
\sup_{w \in \p {\bf d}_\ve} 
   |K_\Omega(p - \ve n_p, w)| = +\infty \, .
$$
This gives a contradiction, and the result is established.
\endpf

\bigskip  

\section{Concluding Remarks}

The worm domains are assuming an ever more prominent role in
the function theory of several complex variables.  Originally
created for the study of elementary facts about the geometry
of pseudoconvex domains, they are now playing an ever-more-prominent
role in the hard analytic questions of these domains.  It is becoming
clear that understand the worms will help us to understand pseudoconvex
domains in general.  Particularly, it is now apparent that the right
approach to the function theory of a domain is best formulated in the
language of the invariant geometry of that domain.  The subtleties of
the worms will bring that geometry to the fore, and help us to push
this program forward.

\end{document}